\documentclass[11pt,notitlepage,twoside]{article}
\usepackage{latexsym}
\usepackage{amsmath}
\usepackage{amsfonts}
\usepackage{amssymb}
\renewcommand{\a }{\alpha }
\renewcommand{\b }{\beta }
\renewcommand{\d}{\delta }

\newcommand{\D }{\Delta }

\newcommand{\e }{\varepsilon }
\newcommand{\g }{\gamma}

\renewcommand{\l }{\lambda }
\renewcommand{\L }{\Lambda }

\newcommand{\n }{\nabla }
\newcommand{\var }{\varphi }

\newcommand{\Sig }{\Sigma}

\renewcommand{\th }{\theta }

\renewcommand{\O }{\Omega }

\newcommand{\ov}{\overline}
\newcommand{\wtilde }{\widetilde}

\newcommand{\be}{\begin{equation}}
\newcommand{\ee}{\end{equation}}
\newenvironment{pf}{\noindent{\bf Proof.}\enspace}{
\hfill$\Box$\medskip}
\newenvironment{pfn}[1]{\noindent{\bf Proof of {#1}\enspace}}{
\hfill$\Box$\medskip}
\newcommand{\R}{\mathbb{R}}
 
\newcommand{\N}{\mathbb{N}}

\newtheorem{thm}{Theorem}[section]
\newtheorem{pro}[thm]{Proposition}
\newtheorem{lem}[thm]{Lemma}
\newtheorem{rem}[thm]{Remark}

\numberwithin{equation}{section}

\author{Mohamed BEN AYED$^a$\thanks{E-mail: Mohamed.Benayed@fss.rnu.tn} \& Khalil EL MEHDI$^{b,c}$\thanks{Corresponding author. E-mail : \texttt{elmehdik@ictp.trieste.it}, \texttt{khalil@univ-nkc.mr}} \\
{\footnotesize
a : D{\'e}partement de Math{\'e}matiques, Facult{\'e} des Sciences de Sfax, Route Soukra, Sfax, Tunisia.}\\
{\footnotesize
b :  Facult\'e des Sciences et Techniques, Universit\'e de Nouakchott, BP 5026, Nouakchott, Mauritania.}\\
{\footnotesize
 c : The Abdus Salam ICTP, Mathematics Section, Strada Costiera 11, 34014 Trieste, Italy.}
}

\title { \Large \textbf{The  Paneitz Curvature Problem on
     Lower Dimensional Spheres }}

\begin{document}

\date{ }

\maketitle

{\footnotesize

\noindent
{\bf Abstract.}
In this paper we prescribe a fourth order
conformal invariant (the
Paneitz curvature) on the $n$-spheres, with $n\in\{5,6\}$. Using dynamical and
topological methods involving the study of critical points at infinity
of the associated variational problem, we prove some existence
results.
}

\section{Introduction and the Main Results}
\mbox{}
In \cite{P}, Paneitz introduced a conformally fourth order operator defined on $4$-manifolds. In \cite{Br1}, Branson generalized the definition to $n$-dimensional Riemannian manifolds, $n\geq 5$. Given a smooth compact Riemannian $n$-manifold $(M,g)$, $n\geq 5$, let $P^n_g$ be the operator defined by
$$
P^n_g u= \D ^2_gu - div_g(a_n S_g g + b_nRic_g ) du + \frac{n-4}{2} Q^n_g u,
$$
where
$$
a_n=\frac{(n-2)^2 +4}{2(n-1)(n-2)}, \qquad b_n =\frac{-4}{n-2}
$$
$$
Q^n_g = - \frac{1}{2(n-1)} \D _g S_g + \frac{n^3 -4n^2 + 16n-16}{8(n-1)^2(n-2)^2}S^2_g - \frac{2}{(n-2)^2} |Ric_g|^2
$$
and where $S_g$ denotes the scalar curvature of $(M,g)$ and
$Ric_g$ denotes the Ricci curvature of $(M,g)$.\\
Such a $Q^n$ is a fourth order invariant and we call it the
Paneitz curvature. For more details about the properties of the
Paneitz operator, see for example \cite{Br1}, \cite{Br2},
\cite{BCY}, \cite{C}, \cite{CGY1}, \cite{CQY1}, \cite{CQY2},
\cite{CY1}, \cite{DHL}, \cite{DMO1}, \cite{DMO2}, 
\cite{F}, \cite{G} and the references therein.\\
If $\tilde{g}=u^{4/(n-4)}g$ is a metric  conformal to $g$,
where $u$ is a smooth positive function,  then for all
$\var \in C^\infty (M)$ we have
$$
P^n_g(u\var ) = u^{(n+4)/(n-4)} P^n_{\tilde{g}}(\var ).
$$
Taking $\var \equiv 1$, we then have
\begin{eqnarray}\label{e:1}
P^n_g(u) =\frac{n-4}{2} Q^n_{\tilde{g}} u^{(n+4)/(n-4)}
\end{eqnarray}
In view of equation \eqref{e:1}, a natural question
is whether it is possible to prescribe the Paneitz
curvature, that is: given a function $ f: M \to \R$,
does there exist a metric $\tilde{g}$ conformally
equivalent to $g$ such that $Q^n_{\tilde{g}}=f$ ?
According to equation \eqref{e:1}, the problem is
equivalent to finding a smooth positive solution of
the following equation
\begin{eqnarray}\label{e:2}
P^n_g(u) =\frac{n-4}{2} f u^{(n+4)/(n-4)}, \qquad u > 0\qquad \mbox{ in } M.
\end{eqnarray}
In this paper, we consider the case of standard sphere $S^n$ endowed with its standard metric $g_0$, and in particular the cases $n=5$ and $n=6$. We are thus reduced to find a positive solution $u$ of the problem
\begin{eqnarray}\label{e:3}
 \mathcal{P}u:=\D^2u-c_n\D u+d_n u =K u^\frac{n+4}{n-4}, \qquad u > 0\qquad \mbox{ in } S^n,
\end{eqnarray}
where $c_n=\frac{1}{2}(n^2-2n-4)$ and $d_n=\frac{n-4}{16}n(n^2-4)$ and where $K$ is a given $C^3$ function defined on $S^n$.\\
 More precisely, our aim is to give sufficient conditions on $K$ such that equation \eqref{e:3} possesses a solution. It is easy to see that a necessary condition on $K$ for solving equation \eqref{e:3} is that $K$ has to be positive somewhere. In addition, there are topological obstructions of Kazdan-Warner type to solve \eqref{e:3} (see \cite{DHL} and \cite{WX}) and so a natural question arises : under which conditions on $K$, \eqref{e:3} has a solution. Our aim is to handle such a question, using some topological and dynamical tools of the theory of critical points at infinity, see Bahri \cite{B1}.\\
To state our main results, we need to introduce some notations. Throughout this paper $K$ denotes a positive $C^3$ function on $S^n$ ( $n=5,6$) which has only nondegenerate critical points $y_1$,..., $y_N$ such that $ -\D K(y_i) \neq 0$ for any $i = 1$,...,$N$. Each $y_i$ is assumed to be of Morse index $k_i$. For the sake of simplicity, we assume that
$$
-\D K(y_i) > 0\qquad \mbox{ for }  1\leq i\leq l \qquad \mbox{ and } -\D K(y_i ) < 0 \qquad \mbox{ for } l+1 \leq i\leq N.
$$
For any $s\in\{1,...,l\}$ and for any $s$-tuple
$\tau_s=(i_1,...,i_s)\in \{1,...,l\}^s$ such that
$i_p\neq i_q$ for $p\neq q$, we introduce a matrix
$M(\tau_s)=(m_{pq})_{1\leq p,q \leq s}$ with
\begin{eqnarray}\label{mat}
m_{pp}=\frac{-\D K(y_{i_p})}{K(y_{i_p})^\frac{3}{2}}, \qquad
m_{pq}=-30\frac{G(y_{i_p},y_{i_q})}
{(K(y_{i_p})K(y_{i_q}))^{\frac{1}{4}}}, \quad \mbox{ if } p\ne q,
\end{eqnarray}
where $G$ is a Green function for $\mathcal{P}$ on $S^6$. It is given by $G(x,y)= (1-cosd(x,y))^{-1}$.\\
 Let $Z$ be a pseudogradient of $K$, of
Morse-Smale type (that is, the intersections of the stable and the
unstable manifolds of the critical points of $K$ are transverse.)\\
Set
$$
X=\overline{\cup_{1\leq i\leq l}W_s(y_i)},
$$
where $W_s(y_i)$ is the stable manifold of $y_i$ for $Z$.\\
Now, we are able to state our main results.
\begin{thm}\label{t:2}
Let $n=5$. Assume that the following two assumptions hold :\\
{\bf $(A_1)$}\hskip 0.3cm $X$ is not contractible \\
{\bf $(A_2)$}\hskip 0.3cm $W_s(y_i) \cap W_u(y_j) = \emptyset $ for
any $i\in \{1,...,l\}$ and for any $j\in \{ l+1,...,m\}$.\\
Then \eqref{e:3} has a solution.
\end{thm}
\begin{thm}\label{t:4}
Let $n=6$ and assume that the following assumption holds :\\
  (H) \hskip 0.3cm $\D K(y_i)\D K(y_j) < 900 G(y_i,y_j)^2 K(y_i)
  K(y_j)$ for any $i\ne j$ in $\{1,...,l\}$.\\
 If
$$
1\ne \sum_{s=1}^l(-1)^{k_s},
$$
 where $k_s$ is the Morse index of $K$ at $y_s$, then \eqref{e:3} has a solution.
\end{thm}
\begin{rem}
 The assumption $(H)$ implies that the least
eigenvalue $\rho(y_i,y_j)$ of the matrix \\$M(y_i,y_j)$ defined by
\eqref{mat} is negative for any indices $i\ne j$ in $\{1,...,l\}$.
\end{rem}
\begin{thm}\label{t:5}
Let $n=6$. Assume that the following two assumptions hold :\\
 $(H_1)$\hskip 0.3cm $X$ is not contractible \\
 $(H_2)$\hskip 0.3cm $\D K(y_i)\D K(y_j) < 900 G(y_i,y_j)^2 K(y_i)
  K(y_j)$ for any $i\ne j$ in $\{1,...,l\}$.\\
If\\
$W_s(y_i) \cap W_u(y_j) = \emptyset $ for any $i\in \{1,...,l\}$ and for any $j\in \{ l+1,...,m\}$,\\
then \eqref{e:3} has a solution.
\end{thm}

Our approach extends the topological and dynamical methods
developed by Bahri \cite{B2}, Bahri-Coron \cite{BC1} and Ben Ayed
et al \cite{BCCH}, to the framework of such higher order
equations. To perform such an extension, a fine analysis of the
gradient flow of the Euler Lagrange Functional is needed. It turns
out that such a gradient flow satisfies the Palais smale condition
on its decreasing flow lines far from a finite number of isolated
blow up. Then we construct a special pseudogradient near such
"singularities" and perform a Morse reduction. Such a fine
analysis of these "singularities" has its own interest, and plays
a central role in the derivation of further existence results to
be
published in forthcoming papers.\\
Another main issue in our approach is to prove the positivity of the
critical point obtained by our process. It is known that in the
framework of such a higher order equation, such an issue is far from been
trivial in general (see \cite{DMO2} for example), and the way we handle it here is very simple compared with the
literature.

Besides the above results, we point out that our method enables us to reprove some existence results, obtained recently by Djadli, Malchiodi and Ould Ahmedou \cite{DMO2}, namely:
\begin{thm}\label{t:1}
Assume that $n=5$. If
$$\sum_{1\leq i\leq l}(-1)^{k_i}\ne -1,$$
where $k_i$ is the Morse index of $K$ at $y_i$, then \eqref{e:3} has a solution.\end{thm}
\begin{thm}\label{t:3}
Let $n=6$ and assume that for any $s\in \{1,...,l\}$ and for any  $\tau_s$,
$M(\tau_s)$ is nondegenerate. If
$$
1\ne \sum_{s=1}^l\sum_{\tau_s=(i_1,...,i_s)/\rho(\tau_s)>0}(-1)^{7s-1-
\sum_{j=1}^sk_{i_j}},$$
then \eqref{e:3} has a solution, where $\rho (\tau _s)$
denotes the least eigenvalue of $M(\tau _s)$.
\end{thm}

We organize  our paper as follows. In section 2, we set up
the variational structure and recall some preliminaries.
In section 3, we perform an expansion of the Euler functional
associated to \eqref{e:3} and its gradient near critical
points at infinity, then in section 4, we give the
characterization of the critical points at infinity.
In section 5, we provide the proofs of our results.
The proofs require some technical results which, for
the convenience of the reader, are given in the appendix.
\section{Preliminary Tools }
\mbox{}
In this section we recall the functional setting and
the variational problem and its main features.\\
For $K\equiv1$, the solutions of \eqref{e:3} are the family
$\wtilde{\d}_{(a,\l)}$ defined by
$$
\wtilde{\d}_{(a,\l)}(x)=\b_n \frac{1}{2^{\frac{n-4}{2}}}
\frac{\l^\frac{n-4}{2}}{\bigl(1+\frac{\l^2-1}{2}(1-cos
d(x,a))\bigr)^\frac{n-4}{2}},
$$
 where $a\in S^n$,  $\l >0$ and $\b_n$ is a positive constant.
 After performing a stereographic projection $\pi$ through the
 point $-a$ as pole, the function $\wtilde{\d}_{(a,\l)}$ is
 transformed into
 $$
\d_{(0,\l)}(y)=\b_n\frac{\l^{\frac{n-4}{2}}}{(1+\l^2\mid
 y\mid^2)^{\frac{n-4}{2}}},
$$
 which is a solution of the problem
$$
 \D^2u= u^\frac{n+4}{n-4} ,\, u>0\,\quad \mbox{ on } \quad \R^n \quad \mbox{(see \cite{Li})}.
$$
The space $ H_2^2(S^n)$ is equipped with the norm :
$$
\mid\mid u\mid\mid^2= <u,u>=\int_{S^n}  \mathcal{P}u.u = \int_{S^n} \mid \D u\mid^2 +c_n
\int_{S^n} \mid\n u\mid^2 + d_n\int_{S^n} u^2.
$$
We denote by $\Sig$ the unit sphere of $H_2^2(S^n)$ and we set
$\Sig ^+ =\{u\in \Sig/ u >0\}$.\\
We introduce  the following functional defined on $\Sig$ by
$$
J(u)=\frac{1}{(\int_{S^n}
K\mid u\mid^{\frac{2n}{n-4}})^{\frac{n-4}{n}}}=\frac{\mid\mid
u\mid\mid^2}{(\int_{S^n} K\mid u\mid^{\frac{2n}{n-4}})^{\frac{n-4}{n}}}.
$$
 The positive critical points of $J$, up to a multiplicative
constant, are solutions of \eqref{e:3}. The Palais-Smale condition
fails for $J$ on $\Sig ^+$. This failure can be described using
similar arguments as in \cite{BrC}, \cite{L}, \cite{S}.
\begin{pro}\label{p:21}
Assume that $J$ has no critical point in $\Sig ^+$ and
let $(u_k)$ be a sequence in $\Sig ^+$ such that $J(u_k)$ is bounded
and $\n J(u_k)$ goes to 0. Then there exist an integer $p$ and a
sequence $\e _k$ such that $u_k\in V(p,\e _k)$, where $V(p,\e)$ is
defined by
\begin{align*}
V(p,\e) =\{u\in & \Sig / \exists a_1, ..., a_p \in S^n,
\exists \l _1,..., \l _p > \e^{-1}, \exists \a _1,..., \a _p > 0  \mbox{ with }\\
 & \mid\mid u-\sum_{i=1}^p\a_i\wtilde{\d}_{(a_i, \l_i)}
\mid\mid < \e ; \quad  \mid J(u)^{\frac{n}{n-4}} \a _i
^{\frac{8}{n-4}}K(a_i)-1\mid < \e \, \forall i,  \quad \e _{ij} <
\e \, \forall i\neq j \}.
\end{align*}
Here
$$
\e_{ij}=(\frac{\l_i}{\l_j}+\frac{\l_j}{\l_i}+
\frac{\l_i\l_j}{2}(1-cosd(a_i,a_j)))^{-\frac{n-4}{2}}.
$$
\end {pro}
The following result defines a parametrization of the set $V(p,\e )$.
\begin{pro}\label{p:22}
For any $p\in \N^*$, there exists $\e_p>0$ such that, if
$0<\e<\e_p$ and  $u\in V(p,\e)$, then the following minimization
problem
$$
min \{\mid\mid u-\sum_{i=1}^p\a_i\wtilde{\d}_{(a_i, \l_i)}
\mid\mid,\,  \a_i>0,\, \l_i>0 ,\,  a_i\in S^n\}
$$
has a unique solution $(\a,a,\l)=(\a _1,...,\a _p,a_1,...,a_p,\l
_1,...,\l _p)$ (up to permutation). In particular, we can write $u
\in V(p,\e )$ as follows
$$
u=\sum_{i=1}^p\a_i\wtilde{\d}_{(a_i,\l_i)} +v,
$$
where $v\in H^2_2(S^n)$ such that, for any $i=1,...,p$
\begin{eqnarray}\label{V0}
(V_0):\quad
<v,\varphi _i>=0 \mbox{ for } \varphi \in\{ \wtilde{\d}_{(a_i,\l_i)},\,
\partial \wtilde{\d}_{(a_i,\l_i)}/ \partial\l_i,\, \partial
\wtilde{\d}_{(a_i,\l_i)}/\partial (a_i)_j\}\quad  \forall \,
j=1,...,n,
\end{eqnarray}
for some system of coordinates $(a_i)_1,...,(a_i)_n$ on $S^n$ near
$a_i$.
\end{pro}
The proof of Proposition \ref{p:22} is similar, up to
minor modifications, to the corresponding statements in \cite{BC2} and \cite{B2}.
 \section{Expansion of the Functional and its Gradient }
 \mbox{}
In this section, we perform a useful expansion of the functional
associated to \eqref{e:3} and its gradient near a critical point at infinity.
 \begin{pro}\label{p:31}
For $\e$ small enough and $u=\sum_{i=1}^p\a _i\wtilde{\d}_i+v
\in V(p,\e)$, we have the following expansion
\begin{align*}
J(u)= & \frac{(\sum_{i=1}^p\a _i ^2) S_n^{4/n} }{(\sum_{i=1}^p\a_i
^{\frac{2n}{n-4}}K(a_i))^{\frac{n-4}{n}}} \biggl[1-\frac{n-4}{n}
c_2\sum_{i=1}^p \frac{\a_i ^{\frac{2n}{n-4}}}{\sum_{j=1}^p\a_j
^{\frac{2n}{n-4}}K(a_j)S_n}\frac{4\D K(a_i)}{\l_i ^2}\\
 & +\frac{c_1}{S_n}\sum_{i\ne j}\a
 _i\a_j\e_{ij}\biggl(\frac{1}{\sum_{k=1}^p\a_k ^2
 }-\frac{2\a_i ^{\frac{8}{n-4}}K(a_i)}{\sum_{k=1}^p\a_k
 ^{\frac{2n}{n-4}}K(a_k)}\biggr) -f(v)+ \frac{1}{\sum_{i=1}^p\a_i ^2
 S_n}Q(v,v) \\
& +o\biggl( \sum_{i\ne j}\e_{ij}+\sum \frac{1}{\l_i ^2} + \mid\mid
v\mid\mid^{\min(\frac{2n}{n-4},3)}\biggr) \biggr],
\end{align*}
where $c_1$ and $c_2$ are positive constants (defined in Lemmas
\ref{lem2} and \ref{lem3}), $S_n=\int_{R^n}\d^{2n/(n-4)}$,
$$Q(v,v)=\mid\mid v\mid\mid^2- \frac{n+4}{n-4}\frac{\sum_{i=1}^p\a_i ^2
 }{\sum_{i=1}^p\a_i
^{\frac{2n}{n-4}}K(a_i)}
\int_{S^n}K(\sum_{i=1}^p\a_i\wtilde{\d}_i)^{\frac{8}{n-4}}v^2
$$
and
$$ f(v)=\frac{2}{\sum_{j=1}^p\a_j
 ^{\frac{2n}{n-4}}K(a_j)S_n}
 \int_{S^n}K(\sum_{i=1}^p\a_i\wtilde{\d}_i)^{\frac{n+4}{n-4}}v.
 $$
(Here and in the sequel $\wtilde{\d}_i$ denotes $\wtilde{\d}_{(a_i,\l _i)}$).
\end{pro}
\begin{rem}
According to Proposition \ref{p:31}, we see that there is a
difference between the three cases $n=5$, $n=6$ and the
higher dimensions. In the case $n=5$, the interaction between two masses
dominates the self interaction, while for $n=6$, there is
a balance phenomenon, and for $n\geq 7$, the self interaction
dominates the interaction between two masses.
\end{rem}
\begin{pfn}{\bf Proposition \ref{p:31}}
Let us recall that
$$
J(u)=\frac{\mid\mid u\mid\mid^2}{(\int_{S^n}
Ku^{\frac{2n}{n-4}})^{\frac{n-4}{n}}}.
$$
Using Lemmas \ref{lem1} and \ref{lem2} in the Appendix, we have
\begin{align*}
\mid\mid u\mid\mid^2 & =\sum_{i=1}^p\a _i ^2\mid\mid
\wtilde{\d}_i\mid\mid^2+\sum_{i\ne j}\a _i\a
_j<\wtilde{\d}_i,\wtilde{\d}_j>+\mid\mid v\mid\mid^2\\
 & =  \sum_{i=1}^p\a _i ^2 S_n +\sum_{i\ne j}\a _i\a _j
(c_1\e _{ij}+o(\e _{ij}))
 + \mid\mid v\mid\mid^2\\
 & = (\sum_{i=1}^p\a _i ^2 S_n )\biggl(1+c_1\sum_{i\ne j}\frac{\a
 _i\a_j}{\sum_{k=1}^p\a_k ^2 S_n }\e_{ij}+\frac{1}{\sum_{i=1}^p\a_i ^2
 S_n}\mid\mid
 v\mid\mid^2+o(\sum_{i\ne j}\e_{ij})\biggr).
\end{align*}
Furthermore, we have
\begin{align}\label{d}
  & \int_{S^n} K(\sum_{i=1}^p\a_i\wtilde{\d}_i  +v)^{\frac{2n}{n-4}} =
\int_{S^n}K(\sum_{i=1}^p\a_i\wtilde{\d}_i)^{\frac{2n}{n-4}}+\frac{2n}{n-4}
\int_{S^n}K(\sum_{i=1}^p\a_i\wtilde{\d}_i)^{\frac{n+4}{n-4}}v\\
 & +\frac{n(n+4)}{(n-4)^2}
\int_{S^n}K(\sum_{i=1}^p\a_i\wtilde{\d}_i)^{\frac{8}{n-4}}v^2+O\biggl(\int
(\sum\a_i\wtilde{\d}_i)^{\frac{12-n}{n-4}}
\mbox{inf}^3((\sum\a_i\wtilde{\d}_i),v)+
\mid v\mid_{L^{\frac{2n}{n-4}}}^{\frac{2n}{n-4}}\biggr).\notag
\end{align}
Since the Sobolev embedding of $H_2^2(S^n)$ in
$L^{\frac{2n}{n-4}}$ is continuous, then there exists a constant
$c$ such that $$\int \mid v\mid^{(2n)/(n-4)}\leq c \mid\mid
v\mid\mid^{(2n)/(n-4)}.$$
We also have
\begin{align*}
\int_{S^n}K(\sum_{i=1}^p\a_i\wtilde{\d}_i)^{\frac{2n}{n-4}}  = &
\sum_{i=1}^p\a_i ^{\frac{2n}{n-4}}\int K\wtilde{\d}_i
^{\frac{2n}{n-4}} + \frac{2n}{n-4}\sum_{i\ne j}\a_i
^{\frac{n+4}{n-4}}\a_j\int K \wtilde{\d}_i ^{\frac{n+4}{n-4}}
\wtilde{\d}_j\\
 & + O\biggl(\sum_{i\ne j}\int \wtilde{\d}_i
^{\frac{8}{n-4}}\inf(\wtilde{\d}_i,\wtilde{\d}_j)^2\biggr).
\end{align*}
For $n\geq 8$, we have $8/(n-4)\leq 2$ and  using Lemma
\ref{lem4} we find
 \begin{eqnarray}\label{e:32}
\int \wtilde{\d}_i ^{\frac{8}{n-4}}\inf(\wtilde{\d}_i,
\wtilde{\d}_j)^2\leq \int (\wtilde{\d}_i
\wtilde{\d}_j)^{\frac{n}{n-4}}
=O(\e_{ij}^{\frac{n}{n-4}}log\e_{ij}^{-1}).
\end{eqnarray}
 For $n<8$, we have $8/(n-4)>2$ and  using Lemma
 \ref{lem4} we obtain
 \begin{eqnarray}\label{e:33}
\int \wtilde{\d}_i
^{\frac{8}{n-4}}\inf(\wtilde{\d}_i,\wtilde{\d}_j)^2\leq
\int\wtilde{\d}_i ^{\frac{8}{n-4}}\wtilde{\d}_j^2\leq c
\biggl(\int(\wtilde{\d}_i\wtilde{\d}_j)^{\frac{n}{n-4}}
\biggr)^{\frac{2(n-4)}{n}}=O(\e_{ij}^2(log\e_{ij}^{-1})^{\frac{2(n-4)}{n}}).
\end{eqnarray}
Using Lemmas \ref{lem3}, \ref{lem4}, \eqref{e:32} and
\eqref{e:33}, we derive that
\begin{align*}
\int_{S^n}K(\sum_{i=1}^p\a_i\wtilde{\d}_i)^{\frac{2n}{n-4}} & =
\sum_{i=1}^p\a_i ^{\frac{2n}{n-4}}\biggl(K(a_i)S_n+c_2\frac{4\D
K(a_i)}{\l_i ^2}+O(\frac{1}{\l_i ^3})\biggr)\\
 & + \frac{2n}{n-4}\sum_{i\ne j}\a_i ^{\frac{n+4}{n-4}}\a_j\biggl(
c_1K(a_i)\e_{ij}+o(\e _{ij}+\frac{1}{\l_i ^2})\biggr).
\end{align*}
Then \eqref{d} becomes
\begin{align*}
\int_{S^n}K  (\sum_{i=1}^p & \a_i\wtilde{\d}_i +v)^{\frac{2n}{n-4}} =
\sum_{i=1}^p\a_i ^{\frac{2n}{n-4}}\biggl(K(a_i)S_n+c_2\frac{4\D
K(a_i)}{\l_i ^2}\biggr)\\
 &  + \frac{2n}{n-4}\sum_{i\ne j}\a_i ^{\frac{n+4}{n-4}}\a_j
c_1K(a_i)\e_{ij} +\frac{2n}{n-4}
\int_{S^n}K(\sum_{i=1}^p\a_i\wtilde{\d}_i)^{\frac{n+4}{n-4}}v\\
 &  +\frac{n(n+4)}{(n-4)^2}
\int_{S^n}K(\sum_{i=1}^p\a_i\wtilde{\d}_i)^{\frac{8}{n-4}}v^2
  + O\biggl(\mid\mid v\mid\mid^{\inf(\frac{2n}{n-4},3)}\biggr)
  +o\biggl(\sum_{i\ne j}\frac{1}{\l_i ^2}+
\e_{ij}\biggr)
\end{align*}
Thus our result follows.
\end{pfn}\\
As usual in this type of problem, we first deal with the
$v$-part of $u$. Let us introduce the following set
$$
E=\{v/ v \mbox{ satisfies } (V_0) \mbox{ and } ||v|| <\e \},
$$
where $(V_0)$ is defined in \eqref{V0}.
\begin{pro}\label{p:vb}
For any $u=\sum_{i=1}^p\a_i\wtilde{\d}_i \in V(p,\e)$ given, there
exists a unique $\ov{v}=\ov{v}(a,\a,\l)$ which minimizes $J(u+v)$
with respect to $v\in E$. Moreover, we have the following estimate
\begin{eqnarray*}
\mid\mid \ov{v}\mid\mid\leq c\mid\mid f\mid\mid\leq c \biggl(
\sum_{i=1}^p\frac{\mid\n K(a_i)\mid}{\l_i}+\frac{1}{\l_i
  ^2}+ \sum_{i\ne j} \e_{ij}^{\min(1,\frac{n+4}{2(n-4)})}
(log\e_{ij}^{-1})^{\min(\frac{n-4}{n},\frac{n+4}{2n})} \biggr).
\end{eqnarray*}
\end{pro}
Before we prove this result, we give the following proposition,
whose proof is deferred to the Appendix
\begin{pro}\label{p:Q}
For any $\sum_{i=1}^p\a_i\wtilde{\d}_i\in V(p,\e)$ given,
$Q(v,v)$ is a quadratic positive form in the space $E$.
\end{pro}
\begin{pfn}{\bf Proposition \ref{p:vb}}
On  one hand, using Proposition \ref{p:Q}, we derive $||v||
< c||f||$, with $c>0$. On the other hand, we have
 $$ f(v)=2(\sum_{j=1}^p\a_j
 ^{\frac{2n}{n-4}}K(a_j)S_n)^{-1}
 \int_{S^n}K(\sum_{i=1}^p\a_i\wtilde{\d}_i)^{\frac{n+4}{n-4}}v.
 $$
 Observe that
  \begin{align*}
  \int_{S^n} & K(\sum_{i=1}^p\a_i\wtilde{\d}_i)^{\frac{n+4}{n-4}}v=
  \sum_{i=1}^p\a_i ^{\frac{n+4}{n-4}}\int
  K\wtilde{\d}_i ^{\frac{n+4}{n-4}}v+O\biggl(\sum_{i\ne j}\int \wtilde{\d}_i
  ^{\frac{8}{n-4}}\inf (\wtilde{\d}_i,\wtilde{\d}_j)\mid v\mid\biggr)\\
 & =O\biggl(\sum_{i=1}^p \mid \n K(a_i)\mid\int\mid x-a_i\mid\wtilde{\d}_i
 ^{\frac{n+4}{n-4}}\mid v\mid+\frac{\mid\mid v\mid\mid}{\l_i ^2}
+\sum_{i\ne j}\int \wtilde{\d}_i
  ^{\frac{8}{n-4}}\inf (\wtilde{\d}_i,\wtilde{\d}_j)\mid
  v\mid\biggr)\\
  & \leq c\mid\mid v\mid\mid\biggl(\sum_{i=1}^p\frac{\mid \n K(a_i)\mid}{\l_i}+\frac{1}{\l_i
  ^2}+\sum_{i\ne j} \e_{ij}^{\min(1,\frac{n+4}{2(n-4)})}
  (log\e_{ij}^{-1})^{\min(\frac{n-4}{n},\frac{n+4}{2n})}\biggr).
  \end{align*}
 Thus the result follows.
 \end{pfn}
 \begin{pro}\label{p:34}
For any $u=\sum_{i=1}^p\a_i\wtilde{\d}_i \in V(p,\e)$, we have the following expansion
\begin{align*}
<\n J(u), \l_i \frac{\partial \wtilde{\d}_i}{\partial \l_i}>
  & =2J(u)\biggl[-c_1\sum_{j\ne i} \a_j \l_i \frac{\partial
\e_{ij}}{\partial
\l_i} +\frac{n-4}{n}c_2\a_i ^{\frac{n+4}{n-4}}J(u)^{\frac{n}{n-4}}\frac{4\D
K(a_i)}{\l_i ^2}\\
&  +o\biggl( \sum \frac{1}{\l_k ^2}+\sum \e_{kj}\biggr)\biggr].
\end{align*}
\end {pro}
\begin{pf}
 We have
$$
\n J(u)=2J(u)\biggl[u-J(u)^{\frac{n}{n-4}}\mathcal{P}^{-1}
(Ku^{\frac{n+4}{n-4}}) \biggr].
$$
Thus
$$<\n J(u),\l_i \frac{\partial \wtilde{\d}_i}{\partial \l_i}> =
2J(u)\biggl[\sum_{j=1}^p\a_j<\wtilde{\d}_j,\l_i \frac{\partial
\wtilde{\d}_i}{\partial \l_i}> -J(u)^{\frac{n}{n-4}}\int
K(\sum_{j=1}^p\a_j\wtilde{\d}_j)^{\frac{n+4}{n-4}} \l_i
\frac{\partial \wtilde{\d}_i}{\partial \l_i}\biggr].
$$
Observe that
\begin{align}\label{01}
\int K(\sum_{j=1}^p\a_j\wtilde{\d}_j)^{\frac{n+4}{n-4}} & \l_i
\frac{\partial \wtilde{\d}_i}{\partial \l_i}  = \sum_{j=1}^p
\a_j^{\frac{n+4}{n-4}} \int K\wtilde{\d}_j^{\frac{n+4}{n-4}} \l_i
\frac{\partial \wtilde{\d}_i}{\partial \l_i}+\frac{n+4}{n-4}
\sum_{j\ne i} \int K(\a_i\wtilde{\d}_i)^{\frac{8}{n-4}} \l_i
\frac{\partial \wtilde{\d}_i}{\partial \l_i}(\a_j\wtilde{\d}_j)\notag\\
 & +O\biggl((if\, n\geq 8)\sum_{k\ne j}\int (\wtilde{\d}_j
 \wtilde{\d}_k)^{\frac{n}{n-4}}+(if\, n<8)\sum_{k\ne j}\int \wtilde{\d}_j
^{\frac{8}{n-4}}\wtilde{\d}_k^2\biggr).
\end{align}
Thus using Lemmas \ref{lem1}, \ref{lem3}, \ref{lem4}, and the
fact that $J(u)^{n/(n-4)}\a_i ^{8/(n-4)}K(a_i)=1+o(1)$, for each $i$, the result follows.
\end{pf}
\begin{pro}\label{p:35}
For any $u=\sum_{i=1}^p\a_i\wtilde{\d}_i \in V(p,\e)$, we have
\begin{eqnarray*}
<\n J(u),\frac{1}{\l_i} \frac{\partial \wtilde{\d}_i}{\partial
a_i}>
  = -2c_3 J(u)^{\frac{2n-4}{n-4}}\frac{\n
  K(a_i)}{\l_i}+O(\frac{1}{\l_i ^2}+\sum_{j\ne i}\e_{ij}).
\end{eqnarray*}
\end{pro}
\begin{pf} We have
$$<\n J(u),\frac{1}{\l_i} \frac{\partial \wtilde{\d}_i}{\partial a_i}> =
2J(u)\biggl[\sum_{j=1}^p\a_j<\wtilde{\d}_j,\frac{1}{\l_i}
\frac{\partial \wtilde{\d}_i}{\partial a_i}>
-J(u)^{\frac{n}{n-4}}\int
K(\sum_{j=1}^p\a_j\wtilde{\d}_j)^{\frac{n+4}{n-4}} \frac{1}{\l_i}
\frac{\partial \wtilde{\d}_i}{\partial a_i}\biggr].
$$
Furthermore we can obtain \eqref{01} but with $\frac{1}{\l_i}
\frac{\partial \wtilde{\d}_i}{\partial a_i}$ instead of $\l_i
\frac{\partial \wtilde{\d}_i}{\partial \l_i}$. Thus using Lemmas \ref{lem1}
 and \ref{lem3}, the result follows.
\end{pf}
 \section{Characterization of the Critical Points at Infinity }
 \mbox{}
This section is devoted to the characterization of the critical
points at infinity for lower dimensions ($n=5$ and $n=6$). We recall
that the critical points at infinity are the orbits of the flow that
remain in $V(p,\e (s))$, where $\e (s)$ is a given function such
that $\e (s)$ tends to zero when $s$ tends to $+\infty$ (see \cite{B1}).
 \begin{pro}\label{p:41}
 Let $n=5$, for $p\geq 2$, there exists a pseudogradient $W$
 so that the following holds.\\
\noindent There is a constant $c>0$ independent of
$u=\sum_{i=1}^p\a_i \wtilde{\d}_i\in V(p,\e)$ so that
$$\bigl(-\n J(u+\ov{v}),W+\frac{\partial \ov{v}}{\partial
(\a_i,a_i,\l_i)}(W)\bigr)\geq c\biggl(\sum_{i=1}^p\frac{\mid\n
K(a_i)\mid}{\l_i}+\frac{1}{\l_i ^2}+\sum_{i\ne j}\e_{ij}\biggr).
$$
Furthermore, $\mid W\mid$ is bounded and the $\l_i$'s decrease
along the flow lines.
\end{pro}
\begin{pf}
We order the $\l_i$'s, for the sake of simplicity we can assume that:
$\l_1\leq\l_2\leq...\leq\l_p$. Let $I=\{i/\, \l_i\mid\n K(a_i)\mid \geq 1\}$. Set
$$Z_1=-\sum_{i=2}^p2^i\a_i\l_i\frac{\partial
\wtilde{\d}_i}{\partial \l_i}, \qquad
Z_2=\sum_{i\in I}\frac{1}{\l_i}\frac{\partial
\wtilde{\d}_i}{\partial a_i}\frac{\n K(a_i)}{\mid\n K(a_i)\mid}.$$
Using Propositions \ref{p:34} and \ref{p:35}, we derive that
\begin{eqnarray}\label{e:51}
 <-\n J(u),Z_1> \geq c\sum_{k\ne
r}\e_{kr}+O(\sum_{i=2}^p\frac{1}{\l_i
^2})+o(\frac{1}{\l_1 ^2}).
\end{eqnarray}
\begin{eqnarray}\label{e:52}
<-\n J(u),Z_2> \geq c\sum_{i\in I} \frac{\mid\n
K(a_i)\mid}{\l_i}+O(\sum_{k\ne
r}\e_{kr})+O(\sum_{i\in I}\frac{1}{\l_i ^2}).
\end{eqnarray}
Let $\mu >0$ such that,  for any  critical point $y$ of $K$,
if $d(a,y)\leq 2\mu$ then $\mid\D K(a)\mid > c>0$. Two cases may occur.\\
\noindent
 {\bf Case 1}\qquad
$\l_2\leq \l_1^2$ or $d(a_1,y)>\mu$ for any critical point $y$.\\
\noindent
 In this case, we set $W_1=MZ_1+Z_2$ where $M$ is a large
constant. Observe that in the case where $d(a_1,y)>\mu$, we can
appear $1/\l_1$ on the lower bound of \eqref{e:52} and therefore
all the $1/\l_i$'s. Combining
 \eqref{e:51} and \eqref{e:52}, we derive
 \begin{eqnarray}\label{e:53}
 <-\n J(u),W_1>\geq c\biggl(\sum_{i=1}^p\frac{\mid \n
K(a_i)\mid}{\l_i}+\frac{1}{\l_i ^2}+\sum_{i\ne j}\e_{ij}\biggr).
\end{eqnarray}
In the other case, that is, $\l_2\leq \l_1^2$, we can  easily
prove that $\frac{1}{\l_1^2}=o(\e_{12})$ (since we have
$(\l_1\l_2)^{1/2}d( a_1,a_2) \leq c\l_1^{3/2}=o(\l_1^2)$ and
$(\l_2/\l_1)^{1/2}\leq \l_1^{1/2}=o(\l_1^2)$). Therefore we can
also
obtain \eqref{e:53} in this case.\\
 {\bf Case 2}\qquad
  $\l_2\geq \l_1^2$ and $d(a_1,y)\leq 2\mu$ for
a critical point $y$.\\ \noindent We set $Z_3=\mbox{sign}(-\D
K(y))\l_1\frac{\partial \wtilde{\d}_1}{\partial \l_1}$, that is, we
increase $\l_1$ if $-\D K(y)>0$ otherwise we decrease it.
We define $W_2=MZ_1+Z_3+mZ_2$, where $M$ is a large
constant and $m$ is a small constant. Observe that
\begin{align*}
<-\n J(u),Z_3> &\geq \frac{c}{\l_1^2}+O(\sum_{j\ne 1}\e_{1j})\\
 <-\n J(u),W_2> & \geq cM\sum_{k\ne
r}\e_{kr}+o(\frac{1}{\l_1 ^2})+
\frac{c}{\l_1^2}+O(\sum_{j\ne 1}\e_{1j})+ m\sum_{i\in I}
\frac{\mid\n K(a_i)\mid}{\l_i}\\
  & +O(m\sum_{k\ne r}\e_{kr})+O(\frac{m}{\l_1^2})
  \geq c\biggl(\sum_{i=1}^p\frac{|\n
K(a_i)|}{\l_i}+\frac{1}{\l_i ^2}+\sum_{i\ne j}\e_{ij}\biggr).
 \end{align*}
(since $M$ is large and $m$ is small). The pseudogradient $W$
will be
built as a convex combination of $W_1$ and $W_2$. \\
Arguing as in Appendix B of \cite{BCCH}, we easily derive that
\begin{align}\label{+v}
<-\n J(u+\ov{v}), & W+\frac{\partial \ov{v}}{\partial
(\a_i,a_i,\l_i)}(W)>  \notag\\
 & \geq <-\n J(u),W>
+o\biggl(\sum_{i=1}^p\frac{|\n K(a_i)|}{\l_i}+\frac{1}{\l_i
^2}+\sum_{i\ne j}\e_{ij}\biggr)
\end{align}
and therefore the proposition follows under \eqref{+v}.
\end{pf}
\begin{pro}\label{p:42}
For $n= 5$, there exists a pseudogradient $W$ so that the
following holds. \noindent There is a constant $c>0$ independent
of $u=\a\wtilde{\d}_{(a,\l)}
\in V(1,\e)$ such that\\
1) $-<\n J(u),W>\geq c(\frac{\mid\n K(a)\mid}{\l}+\frac{1}{\l^2})$\\
2) $-<\n J(u+\ov{v}),W+\frac{\partial \ov{v}}{\partial
(\a,a,\l)}(W)> \geq c(\frac{\mid\n K(a)\mid}{\l}+\frac{1}{\l^2})$\\
3) $W$ is bounded\\
4) the only region where $\l$ increases along the flow lines of $W$ is the region
where $a$ is near a critical point $y$ of $K$ with $-\D K(y)>0$.
\end{pro}
\begin{pf}
Let $\mu > 0$ such that, for any critical point $y$ of $K$,
if $d(x,y)\leq 2\mu$ then\\ $\mid\D K(x)\mid > c >0$. Two cases may occur.\\
\noindent {\bf Case 1} $d(a,y)>\mu$ for any critical point
$y$. In this case we have $\mid\n K(a)\mid > c >0$. Set
$$Z_1=\frac{1}{\l}\frac{\partial \wtilde{\d}}{\partial
a}\frac{\n K(a)}{\mid \n K(a)\mid}.$$
 From Proposition \ref{p:35}, we have
$$-<\n J(u),Z_1>\geq c\frac{\mid\n K(a)\mid}{\l}+O(\frac{1}{\l^2})\geq
c\frac{\mid\n K(a)\mid}{\l}+\frac{c}{\l^2}.$$
 \noindent
  {\bf Case 2}
  $d(a,y)\leq 2\mu$ where $y$ is a
critical point of $K$. Set
$$
Z_2=\mbox{Sign}(-\D K(y)) \l\frac{\partial \wtilde{\d}}{\partial \l}+
m\var(\l\mid\n K(a)\mid)Z_1,
$$
where $m$ is a small constant and $\var$ is a $C^\infty$ function
which satisfies $\var(t)=1$ if $t\geq 2$ and $\var(t)=0$ if $t\leq
1$. Using Propositions \ref{p:34} and \ref{p:35}, we derive that
$$
<-\n J(u),Z_2>\geq \frac{c}{\l^2}
+cm (\frac{\mid\n K(a)\mid}{\l}+O(\frac{1}{\l^2})) \geq c\frac{\mid\n
K(a)\mid}{\l}+\frac{c}{\l^2}.
$$
Hence $W$ will be built as a convex combination of $Z_1$ and $Z_2$.
 Thus the proof of claim 1) is completed. Claims 3) and 4) can be
 easily derived  from the definition of $W$.
Regarding the estimate 2), it can be obtained, arguing as in \cite{B2}
and \cite{BCCH}, using Claim 1).
\end{pf}
\begin{pro}\label{p:43}
Let $n=5$. Assume that $J$ has no critical point in $\Sig^+$.
Then the only critical points at infinity of $J$ correspond to
$\tilde\d_{(y,\infty)}$, where $y$ is a critical point of $K$
with $-\D K(y) > 0$. Moreover, such a critical point at infinity
has a Morse index equal to $5-index(K,y)$.
\end{pro}
\begin{pf}
Using Proposition \ref{p:21}, we derive that $|\n J| \geq c$ outside
of $\cup _{p\geq 1}V(p,\e )$, where $c$ is a positive constant
which depends on $\e$. From Proposition \ref{p:41}, we easily
deduce the fact that there is no critical point at infinity in
$V(p,\e )$ for $p\geq 2$. It only remains to see what happens
in $V(1,\e )$. From Proposition \ref{p:42}, we know that the only
region where $\l$ increases along the pseudogradient $W$,
defined in Proposition \ref{p:42}, is the region where $a$ is
near a critical point $y$ of $K$ with $-\D K(y) > 0$. Arguing as
in \cite{B2} and \cite{BCCH}, we can easily deduce from Proposition
\ref{p:42} the following normal form :\\
If $a$ is near a critical point $y$ of $K$ with $-\D K(y) >0$,
we can find a change of variable $(a,\l ) \longrightarrow (\tilde{a},
\tilde\l )$ such that
$$
J(\d_{(a,\l)}+\bar{v})=\Psi (\wtilde{a},\wtilde{\l}):=
\frac{S_5^{4/5}}{K(\wtilde{a})^{1/5}}\biggl(1-\frac{(c-\eta)\D K(y)}
{\wtilde{\l}^2 K(y)}\biggr),
$$
where $c$ is a positive constant and $\eta$ is a small positive constant.\\
This yields a split of variables $a$ and $\l$, thus it is easy to see
that if $\wtilde{a}=y$, only $\wtilde{\l}$ can move. To decrease
the functional $J$, we have to increase $\wtilde{\l}$, thus we
obtain a critical point at infinity only in this case and our result follows.
\end{pf}\\
Next, we are going to study the case when $n=6$. For this purpose,
we divide the set $V(p,\e)$ into five sets.
\begin{align*}
V_1(p,\e,\eta) & =\{u/a_i\in B(y_{j_i},\eta),\, j_i\ne
j_k\,\mbox{ for } \, i\ne k\, \mbox{ and  for }\, \tau=(j_1,...,j_p), \,
\rho(\tau)>0\}\\
V_2(p,\e,\eta) & =\{u/  a_i\in B(y_{j_i},\eta) ,\, j_i\ne
j_k\, \mbox{ for } \, i\ne k,\, -\D K(y_{j_i}) >0, \, \,
\rho(j_1,...,j_p)<0\}\\
V_3(p,\e,\eta) & =\{u/a_i\in B(y_{j_i},\eta),\, j_i\ne j_k\,
\mbox{ for } \, i\ne k,\, \exists \, \, j_1,...,j_r\, s.
t.\, -\D K(y_{j_i})<0\}\\
V_4(p,\e,\eta) & =\{u/a_i\in B(y_{j_i},\eta),\, \exists
\, i\ne k\, \mbox{ such  that }\,  j_i= j_k \}\\
V_5(p,\e,\eta) & =\{u/\exists\,  i_1,...,i_q\, \mbox{ such
that }\, \mid a_{i_j}-y\mid >\eta/2 \, \mbox{ for  all critical
points }\, y\}
\end{align*}
where $\eta$ is a positive constant such that $\eta < (1/4)\inf _{i\ne j}
d(y_i,y_j)$ and for each $i$, if $d(x,y_i)\leq \eta$ then we have $|\D K(x)|> C >0$.\\
 We then have  the following crucial result.
 \begin{pro}\label{p:44}
Let $n=6$, for $p\geq 1$, there exists a pseudogradient $W$
so that the following holds.\\
\noindent There is a constant $c>0$ independent of
$u=\sum_{i=1}^p\a_i \wtilde{\d}_i\in V(p,\e)$ so that \\
(i) $\displaystyle{\bigl(-\n J(u),W\bigr)\geq
c\biggl(\sum_{i=1}^p\frac{\mid\n K(a_i)\mid}{\l_i}+\frac{1}{\l_i
^2}+\sum_{i\ne j}\e_{ij}\biggr).} $\\
 \noindent
 (ii) $\displaystyle{\bigl(-\n
J(u+\ov{v}),W+\frac{\partial \ov{v}}{\partial
(\a_i,a_i,\l_i)}(W)\bigr)\geq c\biggl(\sum_{i=1}^p\frac{\mid \n
K(a_i)\mid}{\l_i}+\frac{1}{\l_i
^2}+\sum_{i\ne j}\e_{ij}\biggr).} $\\
 \noindent
 (iii) $\mid W\mid$ is bounded. Furthermore, when $u\in
\cup_{j=3,4,5} V_j$, we have $d\l_i(W)\leq 0$. When $u\in V_1\cup
V_2$, $\mid d\l_i(W)\mid \leq c\l_i$ for any $i$. Moreover, the
only case where the maximum of the $\l_i$'s is not bounded is when $u\in V_1$.
\end{pro}
\begin{pf}
We start by proving Claim (i). By the assumption, for any critical
point $y$, we have $\D K(y)\ne 0$. Thus we can choose $\eta>0$
such that for any $x\in B(y,\eta)$, we have $\mid\D K(x)\mid>c >0$.
We will define the pseudogradient depending on the sets $V_i$ to
which $u$ belongs. \\
\noindent
 First, we consider the case of  $u=\sum_{i=1}^p\a_i \wtilde{\d}_i\in V_1(p,\e,\eta)$, we
 have for any $i\ne j$, $\mid a_i-a_j\mid >c$ and therefore
 $$\e_{ij}=\frac{2}{\l_i\l_j(1-cos
 d(a_i,a_j))}(1+o(1))=\frac{2G(a_i,a_j)}{\l_i\l_j}(1+o(1)),$$
 where $G(a_i,a_j)=(1-\cos d(a_i,a_j))^{-1}$, it is the Green's function of $\mathcal{P}$. Thus
 $$\l_i\frac{\partial
 \e_{ij}}{\partial\l_i}=-\e_{ij}(1+o(1))=-2\frac{G(a_i,a_j)}{\l_i\l_j}(1+o(1)).
 $$
 Observe that, since $u\in V(p,\e)$, we have $\a_i
 ^4J(u)^3K(a_i)=1+o(1)$. Thus, Proposition \ref{p:34} becomes
\begin{eqnarray}\label{dJ}
 <\n J(u),\a_i\l_i\frac{\partial\wtilde{\d}_i}{\partial
 \l_i}>=2J(u)^{\frac{-1}{2}}\biggl[\frac{c_2}{3}\frac{4\D
 K(a_i)}{K(a_i)^{\frac{3}{2}}\l_i ^2}+\sum_{j\ne
 i}\frac{2c_1 G(a_i,a_j)}{(K(a_i)K(a_j))^{\frac{1}{4}}}\frac{1}{\l_i\l_j}+
 o(\sum_k\frac{1}{\l_k^2})\biggr]
\end{eqnarray}
We define $Z_1$ by $Z_1=\sum_{i=1}^p\a_i\l_i(\partial
\wtilde{\d}_i)/(\partial \l_i)$. Thus, we derive
\begin{eqnarray}\label{Z1}
<-\n J(u),Z_1>=c \, {}^T\L M \L +o(\sum_k\frac{1}{\l_k^2})\geq c
\sum_k\frac{1}{\l_k^2}\geq c\sum_k\frac{1}{\l_k^2}+ c\sum_{i\ne
j}\e_{ij},
\end{eqnarray}
where $M$ is the matrix defined by \eqref{mat} and $\L=^T(1/\l_1,...,1/\l_p)$.\\
\noindent We also define
 $$Z_a=\sum_{i=1}^p\var(\l_i|\n K(a_i)|) \frac{1}{\l_i}\frac{\partial
\wtilde{\d}_i}{\partial a_i}\frac{\n K(a_i)}{\mid \n K(a_i)\mid}\quad
\mbox{ and }\quad W_1=CZ_1+Z_a,
$$
where $C$ is a large constant and where $\var$ is a $C^\infty$
function which satisfies $\var(t)=0$ if $t\leq 1$ and $\var(t)=1$
if $t\geq 2$. Using Proposition \ref{p:35} and \eqref{Z1}, we
derive that
\begin{eqnarray}\label{W1}
<-\n J(u), W_1>\geq c\biggl(\sum_i\frac{\mid \n K(a_i)\mid}{\l_i}
+\sum_k\frac{1}{\l_k^2}+\sum_{i\ne j}\e_{ij}\biggr).
\end{eqnarray}
 \noindent
  Secondly, we study the case of $u=\sum_{i=1}^p\a_i
  \wtilde{\d}_i\in V_2(p,\e,\eta)$. Let  $\rho$ be the least
  eigenvalue of $M$. Then, there exists an
  eigenvector $e$ associated to $\rho$ such that $\mid\mid
  e\mid\mid =1$ with $e_i>0$ for all $i$. Let $\g >0$  such that
  for any $x\in B(e,\g)=\{y\in S^{p-1}/\mid\mid y-e\mid\mid\leq
  \g\}$, we have $^Tx M x < (1/2)\rho$. Two cases may occur.\\
  \noindent
 \underline{Case 1.}  $\mid\L\mid^{-1}\L\in B(e,\g)$. In this
 case, we define $W_2=-CZ_1+Z_a$. As in \eqref{Z1} and
 \eqref{W1}, we derive that
 \begin{eqnarray}\label{W2}
 <-\n J(u), W_2>\geq c\biggl(\sum_i\frac{\mid \n K(a_i)\mid}{\l_i}
+\sum_k\frac{1}{\l_k^2}+\sum_{i\ne j}\e_{ij}\biggr).
\end{eqnarray}
 \noindent
 \underline{Case 2.} $\mid\L\mid^{-1}\L\notin B(e,\g)$. In this case, we define
 $$
Z_2=-\sum_{i=1}^p\mid\L\mid\a_i\l_i ^2\frac{\partial
 \wtilde{\d}_i}{\partial\l_i}\biggl[\frac{\mid\L\mid e_i-\L_i}{\mid\mid
 y(0)\mid\mid}-\frac{y_i(0)}{\mid\mid
 y(0)\mid\mid^3}(y(0),\mid\L\mid e-\L)\biggr],
$$
 where $y(t)=(1-t)\L+t\mid\L\mid e$. Define $\L(t)=y(t)/\mid\mid
 y(t)\mid\mid$. Using Proposition \ref{p:34}, it is easy to derive that
 $$<-\n J(u),Z_2>=-c\mid\L\mid^2\frac{\partial}{\partial
 t}(^T\L(t) M\L(t))+o(\sum\frac{1}{\l_k^2}).$$
 Observe that
 $$^T\L(t) M \L(t) =\rho+\frac{(1-t)^2}{\mid\mid
 y(t)\mid\mid^2}(^T\L M \L-\rho \mid\mid\L\mid\mid^2).$$
 Thus, we derive
 $(\partial)/(\partial t)(^T\L(t) M\L(t))<-c$. Therefore for
 $W_2'=CZ_2+Z_a$, we obtain
  \begin{eqnarray*}
 <-\n J(u), W_2'>\geq c\biggl(\sum_i\frac{\mid \n K(a_i)\mid}{\l_i}
+\sum_k\frac{1}{\l_k^2}+\sum_{i\ne j}\e_{ij}\biggr).
\end{eqnarray*}
 \noindent
 Now, we deal with the case of $u=\sum_{i=1}^p\a_i \wtilde{\d}_i\in V_3(p,\e,\eta)$.
 Without loss of generality, we can assume  that $1,...,q$ are the
 indices which satisfy $-\D K(a_i)<0$.\\ Set
 $I=\{i/\l_i\leq (1/10)\inf_{k=1,...,q}\l_{k}\}$ and let $M_I$ be
 the matrix defined by the points $(a_i)_{i\in I}$ (as in \eqref{mat})
 and $\rho_I$ be the least eigenvalue of $M_I$. Define
 $$Z_3=-\sum_{i=1}^q\a_i\l_i(\partial
\wtilde{\d}_i)/(\partial \l_i)$$
 Then, since $\mid a_i-a_j\mid>c$, using \eqref{dJ}, we obtain
\begin{eqnarray}\label{Z3}
 <-\n J(u), Z_3>\geq c\sum_{k=1}^q\bigl(
\frac{1}{\l_k^2}+\sum_{i\ne
k}\frac{G(a_i,a_k)}{\l_i\l_k}\bigr)\geq c\sum_{k\notin I}\bigl(
\frac{1}{\l_k^2}+\sum_{i\ne k}\e_{ik}\bigr).
\end{eqnarray}
If $I\ne \varnothing$, then the lower bound becomes limited to
those indices such that $k\notin I$. We have to add another
vector field. If the matrix $M_I$ is positive definite, we define
$Z_3'=Z_1(\sum_{i\in I}\a_i \wtilde{\d}_i)$, that means the action
of $Z_1$ but using only the indices of $I$. In the other case,
that is, the matrix $M_I$ is not positive definite, we define
$Z_3'=Z_2(\sum_{i\in I}\a_i \wtilde{\d}_i)$. In both cases, we have
\begin{eqnarray}\label{Z3'}
 <-\n J(u), Z_3'>\geq c\sum_{k\in I}\bigl(
 \frac{1}{\l_k^2}+\sum_{i\ne k,i\in I}\e_{ik}\bigr)-c\sum_{k\in I,i\notin I}\e_{ik}.
\end{eqnarray}
 Now, we define $W_3=CZ_3+Z_3'+mZ_a$ where $C$ is a large constant
 and $m$ is a small constant. Using \eqref{Z3}, \eqref{Z3'} and
 Proposition \ref{p:35}, we derive that
  \begin{eqnarray*}
 <-\n J(u), W_3>\geq c\biggl(\sum_i\frac{\mid \n K(a_i)\mid}{\l_i}
+\sum_k\frac{1}{\l_k^2}+\sum_{i\ne j}\e_{ij}\biggr).
\end{eqnarray*}
 \noindent
 Next, we will study the case of $u=\sum_{i=1}^p\a_i \wtilde{\d}_i\in V_4(p,\e,\eta)$.
Let $B_i=\{j/a_j\in B(y_i,\eta)\}$. In this case, there is at
least one $B_i$ which contains at least two indices. Without loss
of generality, we can assume that $1,...,q$ are the indices such
that the set $B_i$ ($1\leq i\leq q$) contains at least two
indices. We will decrease the $\l_i$'s for $i\in B_i$ with
different speed. For this purpose, let $\chi$ be a smooth cutoff
function such that $\chi\geq 0$, $\chi=0$ if $t\leq \g'$ and
$\chi=1$ if $t\geq 1$, where $\g'$ is a small constant. For $j\in
B_k$, set $\ov{\chi}(\l_j)=\sum_{i\ne j,i\in B_k}\chi(\l_j/\l_i)$.
Define
 $$Z_4=-\sum_{k=1}^q\sum_{j\in
 B_k}\a_j\ov{\chi}(\l_j)\l_j\frac{\partial
 \wtilde{\d}_j}{\partial\l_j}.$$
 Using Proposition \ref{p:34}, we obtain
 \begin{align*}
<-\n J(u),Z_4>=2J(u)\sum_{k=1}^q\sum_{j\in
 B_k}\a_j\ov{\chi}(\l_j)\biggl[ & -c_1\sum_{i\ne
 j}\a_i\l_j\frac{\partial \e_{ij}}{\partial\l_j}+
 \frac{c_2}{3}\a_j^5J(u)^3\frac{4\D K(a_j)}{\l_j^2}\\
 & + o(\sum_r\frac{1}{\l_r^2}+\sum_{i\ne r}\e_{ir})\biggr].
 \end{align*}
 For $j\in B_k$, with $k\leq q$, if $\ov{\chi}(\l_j)\ne 0$, then
 there exists $i\in B_k$ such that $\l_j^{-2}=o(\e_{ij})$ (for
 $\eta$ small enough).\\
  \noindent Furthermore, for $j \in B_k$, if $i\notin B_k$ or
  $i\in B_k$ with $\l_i$ and $\l_j$ are of the same order, that is,
  $\g'<\l_i/\l_j<1/\g'$, then we have  $-\l_r(\partial
  \e_{ij})/(\partial\l_r)=\e_{ij}(1+o(1))$, for $r=i,j$. In the
  case where $i\in B_k$ with (assuming that $\l_i<\l_j$)
  $\l_i/\l_j<\g'$, we have
  $\ov{\chi}(\l_j)-\ov{\chi}(\l_i)\geq 1$. Thus
  $$-\ov{\chi}(\l_j)\l_j(\partial
  \e_{ij})/(\partial\l_j)-\ov{\chi}(\l_i)\l_i(\partial
  \e_{ij})/(\partial\l_i)\geq -\l_j(\partial
  \e_{ij})/(\partial\l_j)=\e_{ij}(1+o(1)).$$
  Thus, we derive that
 \begin{eqnarray}\label{Z4}
  <-\n J(u),Z_4>\geq c\sum_{k=1}^q\sum_{j\in
  B_k,\ov{\chi}(\l_j)\ne 0}\bigl(\frac{1}{\l_j^2}+\sum_{i\ne
  j}\e_{ij}\bigr).
 \end{eqnarray}
  The lower bound does not contain all the indices. We need to add
  some terms. Let
  \begin{eqnarray}\label{i0}
  \l_{i_0}=\inf\{\l_i, \, i=1,...,p\}.
 \end{eqnarray}
Now, we distinguish two subcases.\\
\underline{Subcase 1.} There exists $j$ such that
$\ov{\chi}(\l_j)\ne 0$ and $\l_{i_0}/\l_j>\g'$, then we can
appear on the lower bound $1/\l_{i_0}^2$ and therefore
all the $1/\l_i ^2$ and the $\e_{ik}$. Thus, we can define $W_4^1=CZ_4+Z_a$
 where $C$ is a large constant.\\
\underline{Subcase 2.} For each $j$, we have  $\ov{\chi}(\l_j)= 0$
or $\l_{i_0}/\l_j\leq \g'$. In this case, we define
  $$D=(\{i/\ov{\chi}(\l_i)=0\}\cup(\cup_{k=1}^qB_k)^c)\cap
  \{i\, t.q.\l_i/\l_{i_0}<1/\g'\}.$$
 It is easy to see that $\{i/\ov{\chi}(\l_i)=0\}$ contains at most
 one index from each $B_j$ for $1\leq j\leq q$ and therefore for
 $i,r\in D$ such that $i\ne r$ we have $a_i\in B(y_{j_i},\eta)$
 and $a_r\in B(y_{j_r},\eta)$ with $j_i\ne j_r$. Let
 $$u_1=\sum_{i\in D}\a_i\wtilde{\d}_i.$$
  $u_1$ has to satisfy one
 of the three cases above, that is, $u_1\in V_i(card (D),\e,\eta)$ for
 $i=1,2$ or $3$. Thus, we can apply the associated vector field
 which we will denote $Z_4'$ and we have the estimate
  \begin{eqnarray*}
 <-\n J(u), Z_4'>\geq c\sum_{i\in D}\biggl(\frac{\mid \n K(a_i)\mid}{\l_i}
+\frac{1}{\l_i ^2}+\sum_{i,j\in
D}\e_{ij}\biggr)+O(\sum_{k\in D,r\notin D}\e_{kr})+o(\sum_{i\notin
D}\frac{1}{\l_i ^2}) .
\end{eqnarray*}
Observe that for $k\in D$ and $r\notin D$, we have either $r\in
B^q:=\{i/\ov{\chi}(\l_i)\ne 0\}\cap (\cup_{j=1}^qB_j)$ or $r\in
(B^q)^c$. If $r \in B^q$,
 we have $\e_{kr}$ in the lower bound of \eqref{Z4}. If $r\in
(B^q)^c$, in this case since $r\notin D$ we have
$\l_{i_0}/\l_r<\g'$. Furthermore, we can prove that $a_k$ and $a_r$
are not in the same set $B(y,\eta)$ and therefore $\mid
a_k-a_r\mid>c$. Thus
$$\e_{kr}\leq\frac{c}{\l_k\l_r}\leq
\frac{c\g'}{\l_k\l_{i_0}}=o(\e_{ki_0})$$ ( $\g'$ small). Since
$i_0\in D$ ($i_0$ is defined by \eqref{i0}), then from
$1/\l_{i_0}^2$, we can  appear on the lower bound all the
$1/\l_i ^2$ and $\e_{ir}$ for $i,r\in (B^q)^c$ (since for those
indices we have $\mid a_i-a_r\mid>c$). Thus, we derive that
\begin{eqnarray}\label{Z4'}
<-\n J(u), Z_4'>\geq c\biggl(\sum_{i\in D}\frac{\mid \n
K(a_i)\mid}{\l_i} +\sum_{i=1}^p\frac{1}{\l_i ^2}+\sum_{i,j\in
(B^q)^c}\e_{ij}\biggr)+O(\sum_{k\in D,r\in B^q}\e_{kr}).
\end{eqnarray}
Now, we define $W_4^2=CZ_4+Z_4'+mZ_a$ where $C$ is a large constant
and $m$ is a small constant. We obtain
\begin{eqnarray}\label{W4}
<-\n J(u), W_4^2>\geq c\biggl(\sum_{i=1}^p\frac{\mid \n
K(a_i)\mid}{\l_i} +\sum_{i=1}^p\frac{1}{\l_i ^2}+\sum_{i\ne
j}\e_{ij}\biggr).
\end{eqnarray}
The vector field $W_4$ defined in $V_4(p,\e,\eta)$ will be a convex
combination of $W_4^1$ and $W_4^2$.\\
\noindent
 Finally, we consider the case of $u=\sum_{i=1}^p\a_i \wtilde{\d}_i\in V_5(p,\e,\eta)$.
 We order the $\l_i$'s in an increasing order: $\l_1\leq
 \l_2\leq...\leq \l_p$. Let $i_1$ be such that for any $i<i_1$, we
 have $\mid a_i-y_{j_i}\mid\leq \eta/2$ where $y_{j_i}$ is a critical
 point of $K$ and $\mid a_{i_1}-y\mid>\eta/2$, for any critical point $y$.
 Let us define
 $$u_1=\sum_{i<i_1}\a_i \wtilde{\d}_i$$
 Observe that $u_1$ has to satisfy one of the four cases above
 that is, $u_1\in V_i(i_1-1,\e,\eta)$, for $i=1,2,3$ or $4$. Thus,
 we can apply the associated vector field which we will denote
 $Z_5$ and we have the following estimate
 \begin{align*}
<-\n J(u), Z_5> & \geq c\sum_{i<i_1}\biggl(\frac{\mid \n
K(a_i)\mid}{\l_i} +\frac{1}{\l_i ^2}+\sum_{j<i_1
}\e_{ij}+O(\sum_{j\geq i_1}\e_{ij})\biggr)\\
 & +o(\sum_{i\geq i_1}
\frac{1}{\l_i ^2}+\sum_{k\ne r}\e_{kr}).
\end{align*}
We also define
 \begin{eqnarray*}
 Z_5'=\frac{1}{\l_{i_1}}\frac{\partial\wtilde{\d}_{i_1}}{\partial a_{i_1}}
 \frac{\n K(a_{i_1})}{\mid\n K(a_{i_1})\mid}-C'\sum_{i\geq i_1}2^i\l_i
 \frac{\partial\wtilde{\d}_i}{\partial \l_i},
 \end{eqnarray*}
 where $C'$ is a large constant. Using Propositions \ref{p:34} and \ref{p:35}
 and the fact that $\mid \n K(a_{i_1})\mid >c$, we find
 \begin{align}\label{Z5'}
 <-\n J(u),Z_5'> & \geq \frac{c}{\l_{i_1}}+O(\sum_{i\ne
 i_1}\e_{ii_1})+ cC'\sum_{i\geq i_1}\biggl(\sum_{j\ne
 i}\e_{ij}+O(\frac{1}{\l_i ^2})+o(\sum_{k\ne
 r}\e_{kr}+\frac{1}{\l_k^2})\biggr)\notag\\
  & \geq \sum_{i\geq i_1}(\frac{c}{\l_i}+\sum_{j\ne i}\e_{ij})+ o(\sum_{k\ne
 r}\e_{kr}+\frac{1}{\l_k^2})
 \end{align}
 (since $C'$ is large). Define $W_5=Z_5+CZ_5'$, where $C$ is a
 large constant. We derive that
 \begin{eqnarray}\label{W5}
<-\n J(u), W_5>\geq c\sum_{i=1}^p\biggl(\frac{\mid \n
K(a_i)\mid}{\l_i} +\frac{1}{\l_i ^2}+\sum_{j\ne i }\e_{ij}\biggr).
\end{eqnarray}
Now, we define the pseudogradient $W$ as a convex combination of
$W_i$ for $i=1,...,5$. The construction of $W$ is completed. It
satisfies Claims (i) and (iii) of  Proposition \ref{p:45}.
Regarding (ii), it can be obtained, arguing as in \cite{B2} and
\cite{BCCH},  using the estimate (i).
\end{pf}
\begin{pro}\label{p:45}
Let $n=6$. Assume that $J$ has no critical point in $\Sig ^+$.
Then the only critical points at infinity of $J$ in $\Sig ^+$ correspond to
$$
 \sum _{j=1}^p K(y_{i_j})^{-1/4}\wtilde\d _{(y_{i_j}, \infty)}
 \quad \mbox{ with } \rho(\tau_p)>0 \quad \mbox{ and } \tau_p=(i_1,...,i_p)
$$
where $p\geq 1$ and $\rho (\tau _p)$ denotes the least eigenvalue
of $M(\tau _p)$. Moreover, such a critical point at infinity has a
Morse index equal to $(7p-1-\sum _{j=1}^p \mbox{index }(K,y_{i_j}))$.
\end{pro}
\begin{pf}
>From Proposition \ref{p:21}, we derive that we just need to see
what happens in $V(p,\e )$ for $p\geq 1$. From Proposition
\ref{p:44}, we deduce that the only region where the $\l _i$'s
are not bounded is when each $a_j$ is near critical point
$y_{i_j}$ with $i_j\neq i_k$ for $j\neq k$ and the matrix
$M(\tau _p)$ is  positive definite. In this region, arguing
as in \cite{B2} and \cite{BCCH}, we can find a change of variable
$$
(a_1,...,a_p,\l _1,...,\l _p) \longrightarrow
(\tilde{a}_1,...,\tilde{a}_p,\tilde{\l}_1,...,\tilde{\l}_p):=(\wtilde{a},\wtilde{\l})
$$
such that
$$
J(\sum_{i=1}^p\a _i\tilde \d_{(a_i,\l _i)}+\bar{v})=\Psi (\a
_1,...,\a _p,\wtilde{a},\wtilde{\l}):= \frac{\sum \a
_i^2S_6^{2/3}}{(\sum \a _i^6
K(\wtilde{a}_i))^{1/3}}\biggl(1+(c-\eta)\L ^tM(\tau_p )\L\biggr),
$$
where $c$ is a positive constant, $\eta$ is a small positive
constant and $\L ^t=(1/\tilde{\l}_1,...,1/\tilde{\l}_p)$.\\
Thus, we conclude as in the proof of Proposition \ref{p:43}.\\
It remains to compute the Morse index of such a critical point
at infinity. For this purpose, we observe that $M(\tau_p)$ is
positive definite and the function $\Psi$ admits on the variables
$\a_i$'s an absolute degenerate maximum with one dimensional
nullity space. Then the Morse index of such a critical point at
infinity is equal to $(p-1+\sum_{j=1}^p(6-index(K,y_{i_j})))$.
Thus our result follows.
\end{pf}
 \section{Proofs of Theorems}
 \mbox{}
Let us start by proving the following result adapted from \cite{BCCH}.
\begin{lem}\label{l:51}
For $\eta > 0$ small enough, we define the following neighborhood of $\Sig ^+$
$$
V_\eta (\Sig ^+)=\{u\in \Sig / J(u)^{\frac{2n-4}{n-4}}e^{2J(u)}
|u^-|_{L^{2n/(n-4)}}^{8/(n-4)} < \eta\},
$$
where $u^- = \max (0,-u)$.\\
Then, for $n\geq 5$, $V_\eta (\Sig ^+)$ is invariant under the flow generated by $-\n J$.
\end{lem}
\begin{pf}
Suppose $u_0 \in V_\eta (\Sig ^+)$ and consider
$$
\begin{cases}
\frac{du(s)}{ds}= -\n J(u(s)) = -2J(u)\biggr(u-J(u)^{n/(n-4)}
\mathcal{P} ^{-1}(K|u|^{8/(n-4)}u)\biggr)\\
u(0)=u_0
\end{cases}
$$
Then
$$
e^{2\int_0^s J(u)}u(s)= u_0 +2\int_0^se^{2\int_0^tJ(u)}
J(u)^{\frac{2n-4}{n-4}}\mathcal{P} ^{-1}(K|u|^{8/(n-4)}u)dt.
$$
Therefore
$$
u^-(s) \leq e^{-2\int_0^s J(u)}\biggr( u_0^- +2\int_0^se^{2\int_0^tJ(u)}
J(u)^{\frac{2n-4}{n-4}}\mathcal{P} ^{-1}(K(u^-)^{\frac{n+4}{n-4}})dt\biggr):=
e^{-2\int_0^s J(u)}f(s).$$
Observe that, each solution of $\mathcal{P}v=g$ with $g\in L^{\frac{2n}{n+4}}$,
using a regularity argument, has to satisfy $v\in H_2^2(S^n)$ and then
$v\in L^{\frac{2n}{n-4}}$. Thus $f\in L^{\frac{2n}{n-4}}$.\\
 Setting
$$
F(s)=e^{-\frac{4n}{n-4}\int_0^s J(u)}|f|_{L^{\frac{2n}{n-4}}}^{\frac{2n}{n-4}},
$$
we have $|u^-(s)|_{L^{\frac{2n}{n-4}}}^{\frac{2n}{n-4}} \leq F(s)$.\\
Now, without loss of generality, we can assume that $u_0^-\ne 0$ and we
want to prove that $F$ is a decreasing function. Observe that
\begin{align*}
F'(s)&= -\frac{4n}{n-4} J(u(s))e^{-\frac{4n}{n-4}\int_0^s J(u)}
|f|_{L^{\frac{2n}{n-4}}}^{\frac{2n}{n-4}}+ e^{-\frac{4n}{n-4}
\int_0^s J(u)}\frac{2n}{n-4} \int_{S^n} f'f^{\frac{n+4}{n-4} }\\
&\leq \frac{2n}{n-4} e^{-\frac{4n}{n-4}\int_0^s J(u)}
\left[-2J(u)|u_0^-|_{L^{\frac{2n}{n-4}}}^{\frac{2n}{n-4} } +
\int_{S^n}f'f^{\frac{n+4}{n-4} }\right].
\end{align*}
Notice that
\begin{eqnarray*}
\big|\int_{S^n}f'f^{\frac{n+4}{n-4}}\big|  \leq c\int_{S^n}|f'|
|u_0^-|^{\frac{n+4}{n-4}}+c\int_{S^n}|f'|\biggl(\int_0^s |f'(t)|\biggr)^{\frac{n+4}{n-4}}.
\end{eqnarray*}
But, we have
$$
\int_{S^n}(u_0^-)^{\frac{n+4}{n-4}}\biggl(e^{2\int_0^sJ(u)}
J(u)^{\frac{2n-4}{n-4}}\mathcal{P} ^{-1}(K(u^-)^{\frac{n+4}{n-4}})\biggr)
\leq CJ(u)^{\frac{2n-4}{n-4}} e^{2\int_0^sJ(u)}
|u_0^-|_{L^{{\frac{2n}{n-4}}}}^{\frac{n+4}{n-4}}
|u^-(s)|_{L^{{\frac{2n}{n-4}}}}^{\frac{n+4}{n-4}}.
$$
and we also have
\begin{align*}
\int_{S^n}|f'(s)|(\int_0^s|f'(t))|^{\frac{n+4}{n-4} } & \leq cs^{\frac{8}{n-4}}
\int_{S^n}|f'(s)|\int_0^s|f'(t)|^{\frac{n+4}{n-4}} \\
 & \leq c s^{\frac{8}{n-4}}\biggl( e^{2sJ(u_0)} J(u_0)^{\frac{2n-4}{n-4}}
 \biggr)^{\frac{2n}{n-4}} \int_0^s |u^-(s)|_{L^{{\frac{2n}{n-4}}}}^{\frac{n+4}{n-4}}
 |u^-(t)|_{L^{{\frac{2n}{n-4}}}}^{(\frac{n+4}{n-4})^2}.
\end{align*}
Hence, if $|u^-(s)|_{L^{\frac{2n}{n-4}}} \leq 5|u_0^-|_{L^{\frac{2n}{n-4} }}$,
for $0\leq s\leq 1$, and using the fact that $u_0\in V_\eta(\Sig^+)$, that is,
$J(u_0)^{\frac{2n-4}{n-4}} e^{2J(u_0)}|u_o^-|_{L^{\frac{2n}{n-4}}}^{\frac{8}{n-4}} <\eta $,
and $\eta$ is small enough, then $F'(s) \leq 0$, for $0\leq s\leq 1$. Therefore
$J(u(s))^{\frac{2n-4}{n-4}} e^{2J(u(s))}|u(s)^-|_{L^{\frac{2n}{n-4}}}^{\frac{8}{n-4}} <\eta $,
and our result follows.
\end{pf}\\
Now, we are ready to prove our theorems.\\
\begin{pfn}{\bf Theorem \ref{t:2}}
Arguing by contradiction, we assume that $J$ has no critical point in  $V_\eta (\Sig ^+)$,
where
$$
V_\eta (\Sig^+)=\{u\in\Sig /e^{2J(u)}J(u)^6\mid u^-\mid_{L^{10}}^8<\eta\},
$$
$\eta$ is a small positive constant and $u^-$ denotes the negative
part of $u$,
that is, $u^-=\max (0,-u)$.\\
 It follows from Proposition \ref{p:43} that the only critical
points at infinity of $J$ in $V_\eta (\Sig ^+)$ correspond to
$\tilde\d _{(y,+\infty )}$, where $y$ is a critical point of $K$
with $-\D K(y) > 0$. It follows that $V_\eta(\Sig ^+)$ retracts by
deformation on $X_\infty = \cup _{y_i/-\D
K(y_i)>0}W_u(y_i)_\infty$ (see sections $7$ and $8$ of \cite{BR}),
where $W_u(y_i)_\infty$ is the unstable manifold at infinity of
such a critical point at infinity. Using Assumption $(A_2)$ and
Proposition \ref{p:42}, we see that $X_\infty$ can be parametrized
by $X \times [A,+\infty[$, where $A$ is
a large positive constant. \\
In addition, we have $X_\infty$ is contractible in $V_\eta (\Sig ^+)$ and
$\Sig ^+$ retracts by deformation on $X_\infty$, therefore $X_\infty$ is
contractible leading to the contractibility of $X$, which is in contradiction
with the assumption $(A_1)$ of our theorem. Thus there exists a critical point
of $J$ in $V_\eta(\Sig ^+)$.\\
 Now, it remains to prove that such a critical point is a positive function.
Let us define the function $w^-$ by the solution of the following
problem
$$\mathcal{P} w^-=-K(x) (u^-)^{\frac{n+4}{n-4}}\quad \mbox{on }\, S^n$$
Since $K(x) (u^-)^{\frac{n+4}{n-4}}\in L^{\frac{2n}{n+4}}$, we see $w^-\in H_2^2$.
Furthermore, we have $w^-\leq 0$. Thus we derive
 \begin{align*}
  \int_{S^n} \mathcal{P} w^-.w^- & =\mid\mid w^-\mid\mid^2=\int_{S^n}
 -K(x) (u^-)^{\frac{n+4}{n-4}}w^-\leq C\mid
 w^-\mid_{L^{2n/(n-4)}}\mid u^-\mid_{L^{2n/(n-4)}}^{(n+4)/(n-4)}\\
  & \leq C\mid\mid w^-\mid\mid \, \mid u^-\mid_{L^{2n/(n-4)}}^{(n+4)/(n-4)}.
  \end{align*}
 Thus, either $||w^-||=0$ and therefore $u^-=0$, or $||w^-||\ne 0$ and we derive
 \begin{eqnarray}\label{w-1}
 \mid\mid w^-\mid\mid\leq C\mid
 u^-\mid_{L^{2n/(n-4)}}^{(n+4)/(n-4)}.
\end{eqnarray}
Furthermore, on  one hand we have
\begin{eqnarray}\label{w-2}
 \int_{S^n}u. \mathcal{P} w^-=\int_{S^n}-u
 K(u^-)^{\frac{n+4}{n-4}}=\int_{S^n}K(u^-)^{\frac{2n}{n-4}}\geq
 c_K\int_{S^n}(u^-)^{\frac{2n}{n-4}}\geq c_K \mid
 u^-\mid_{L^{\frac{2n}{n-4}}}^{\frac{2n}{n-4}}
 \end{eqnarray}
 (since $K$ is  bounded from below by a positive constant), and
 on the other hand, we have,
 \begin{align}\label{w-3}
 \int u.\mathcal{P}w^- & =\int w^-.\mathcal{P} u=\int w^- K\mid u\mid^{\frac{8}{n-4}}u\leq
 \int_{u\leq 0}-w^-K(u^-)^{\frac{n+4}{n-4}}\notag \\
  & \leq \int_{S^n}-w^-K(u^-)^{\frac{n+4}{n-4}}=\int_{S^n}
  w^-.\mathcal{P}w^-=\mid\mid w^-\mid\mid^2.
  \end{align}
  Using \eqref{w-1}, \eqref{w-2} and \eqref{w-3}, we obtain
 $$c_K\mid u^-\mid_{L^{2n/(n-4)}}^{2n/(n-4)}\leq \mid\mid
 w^-\mid\mid^2 \leq C \mid
 u^-\mid_{L^{2n/(n-4)}}^{2(n+4)/(n-4)}.$$
 Observe that $2n/(n-4)<2(n+4)/(n-4)$. Thus, either $u^-=0$ or
 $\mid u^-\mid_{L^{2n/(n-4)}}\geq C$ and
this case cannot occur since by the definition of the
neighborhood of $\Sig ^+$ we have this norm is small.
This completes the proof of our result.
 \end{pfn}\\
 \begin{pfn}{\bf Theorem \ref{t:4}}
By Proposition \ref{p:45} and assumption $(H)$, we derive that the
only critical points at infinity of $J$ in  $V_\eta (\Sig ^+)$
correspond to $\wtilde\d _{(y,\infty )}$, where $y$ is a critical
point of $K$ with $-\D K(y) >0$. We order the critical values of
$K$: $K(y_{i_1})\geq K(y_{i_2})\geq ...\geq K(y_{i_l})$ (those
critical points $y_{i_j}$ satisfy $-\D K(y_{i_j})>0$). Let
$c_r=(S_6)^{4/6}(K(y_{i_r}))^{-1/6}$ be the critical value at
infinity. For the sake of simplicity, we can assume that $c_i$'s
are different. Then,  we have
$$
b_1<\min_{\Sig^+}J=c_1<b_2<c_2<b_3<c_3<...<b_l<c_l<b_{l+1}.
$$
Recall that we already built  in Proposition \ref{p:44} a vector
field $W$ defined in $V(p,\e )$ for $p\geq 1$, $\e$ will be chosen
so that $V(p,\e)\subset V_\eta (\Sig ^+)$. Outside $\bigcup
_{p\geq 1}V(p, \e /2)$, we will use $-\n J$ and our global vector
field $Z$ will be built using a convex combination of $W$ and $-\n
J$. Now, according to  Proposition \ref{p:45}, there is no
critical value  above the level $b_{l+1}$. Let  $J_c=\{u\in V_\eta
(\Sig^+)/J(u)<c\}$. Using the vector field $Z$,  we have
$J_{b_{r+1}}$ retracts by deformation on $ J_{b_r} \cup
W_u(y_{i_r})_\infty$, where $W_u(y_{i_r})_\infty$ is the unstable
manifold at infinity (see sections $7$ and $8$ of \cite{BR}).
Then,  denoting by $\chi$ the Euler-Poincar\'e characteristic, we
have
$$
\chi(J_{b_{r+1}})=\chi(J_{b_r})+(-1)^{6-k_r},
$$
where $k_r=index(K,y_{i_r})$. It is easy to see that
$\chi(J_{b_1})=\chi(\varnothing)=0$ and $\chi(V_\eta (\Sig^+))=1$.
Therefore
 \begin{eqnarray}\label{.}
1=\sum_{r=1}^l(-1)^{6-k_r}=\sum_{r=1}^l(-1)^{k_r}.
 \end{eqnarray}
If \eqref{.} is violated, $J$ has a critical point in $V_\eta(\Sig
^+)$. Arguing as in the proof of Theorem \ref{t:2}, we conclude
that this critical point is a positive function and hence our
theorem follows.
\end{pfn}\\
\begin{pfn}{ \bf Theorem \ref{t:5}}
As in the proof of Theorem \ref{t:4}, we derive that the only
critical points at infinity of $J$ in the $V_\eta (\Sig ^+)$
correspond to $\wtilde\d _{(y,\infty )}$, where $y$ is a critical
point of $K$ with $-\D K(y) >0$. Thus, the sequel of the proof of
our theorem is exactly the same as in the proof of Theorem
\ref{t:2}, so we will omit it.
\end{pfn}\\
 \begin{pfn}{Theorem \ref{t:1}}
Arguing by contradiction, we assume that $J$ has no critical point
in $\Sig^+$. Using Proposition \ref{p:43}, the only critical
points at infinity correspond to $\wtilde{\d}(y,\infty)$, where
$y$ is a critical point of $K$ with $-\D K(y)>0$. Such a critical
point at infinity has a Morse index equal to $(5-index(K,y))$.
Using the same arguments as in the proof of Theorem \ref{t:4}, the
result follows.
\end{pfn}\\
 \begin{pfn}{Theorem \ref{t:3}}
 The proof is the same as the proof of Theorem \ref{t:4}. But
 here, the critical points at infinity correspond to
 $$
\sum_{r=1}^pK(y_{i_r})^{-1/4}\wtilde{\d}(y_{i_r},\infty)\qquad \mbox{ with }\quad
 \rho(\tau _p)>0 \mbox{ and } \tau _p= (y_{i_1},...,y_{i_s})),
$$
where $p\geq 1$ and $\rho (\tau _p)$ denote the least eigenvalue
of $M(\tau _p)$. Such a critical point at infinity has an index
equal to $\sum_{k=1}^p(6-k_{i_r})+(p-1)=7p-1-\sum_{r=1}^p
\mbox{index}( K, y_{i_r})$. Using  the  same argument as the proof
of Theorem \ref{t:4}, the result follows.
 \end{pfn}
 \section{Appendix }
\subsection{The Coercivity of the Quadratic Form}
\mbox{}
In this appendix, we give the proof of Proposition \ref{p:Q}, adapted from \cite{B1}.
\begin{pro}\label{Q}
For any $u=\sum_{i=1}^p\a_i\wtilde{\d}_i\in V(p,\e)$ given,
$Q(v,v)$ is a quadratic positive form in the space $E=\{v/\, v\,
\mbox{ satisfies }\, (V_0)\}$.
\end{pro}
\begin{pf}
Using a stereographic projection, we need to prove the proposition
on $R^n$ with the bilaplacian.\\
\noindent Let us define the sets, for $i=1,...,p$
$$\O_i=\{x\in R^n/\mid x-x_i\mid<\frac{1}{8\l_i}\min
\e_{ij}^{\frac{-1}{n-4}} \mbox{ and } \mid x-x_j\mid>\frac{1}{8\l_j}\min
\e_{ij}^{\frac{-1}{n-4}} \mbox{ for } \l_j\, s.t.\, \l_j\geq
\l_i\}$$ By construction $\O_i\cap\O_j=\varnothing$ for $i\ne j$.
Now, we define
$$ H=\{u\in L^{\frac{2n}{n-4}}(R^n)/\D u\in L^2(R^n)\},$$
$H$ is the completion of $C_c^\infty (\R^n)$ with respect to the norm $\int_{\R^n}|\D u|^2$.\\
 For $ \varphi$ belongs to $H$, we introduce
the projection $Q_i$ by:
$\varphi_i=Q_i\varphi$ satisfies
$$ \Delta ^2\varphi_i=\Delta
^2\varphi\, \mbox{in }\, \Omega_i\quad , \quad
\Delta\varphi_i=\varphi_i=0\, \mbox{on }\,
\partial\Omega_i.$$
 Let us define also $E_i ^-=span<\d_i,(\partial\d_i)/(\partial\l_i),
(\partial\d_i)/(\partial a_i)>$ and $E_i ^+=(E_i ^-)^\bot$ (the
orthogonal being taken in the sense of the scalar product $\int
\D\psi\D\var$).\\
Next, we will use the following lemmas which we will prove in the end.
\begin{lem}\label{53}
If, for $i\ne j$, $\e_{ij}$'s are small enough, then
$$\int_{R^n}\d_i ^{\frac{n+4}{n-4}}\mid\var-\var_i\mid\leq
c(\sum\e_{ij}^{1/2})(\int_{R^n}\mid\D\var\mid^2)^{1/2}.$$
\end{lem}
\begin{lem}\label{54}
There exists $\a_1>0$ s.t. for any $\var\in E_i ^+$, we have
$$\int_{R^n}\mid
\D\var\mid^2-\frac{n+4}{n-4}\int_{R^n}\d_i
^{\frac{8}{n-4}}\var^2\geq \a_1\int_{R^n}\mid\D\var\mid^2.$$
\end{lem}
\begin{lem}\label{55}
For $v\in H$ satisfying $(V_0)$ and $v_i=Q_i v$, we write $v_i=v_i
^-+v_i ^+$, where $v_i ^-\in E_i ^-$ and $v_i ^+\in E_i ^+$. Then,
we have
$$\int_{\O_i}\mid\D v_i ^-\mid^2=o(\int_{R^n}\mid\D v\mid^2).$$
\end{lem}
Using those Lemmas, we are able to give the proof of the above
proposition. Indeed:\\
 Let $v$ satisf $(V_0)$, we denote $v_i=Q_i v$ for each
$i=1,...,p$. We can assume that $v_i$ is defined on $R^n$ by taking
$v_i=0$ on $\O_i ^c$. We split $v_i$
into two parts: $v_i=v_i ^-+v_i ^+$ where $v_i ^-\in E_i ^-$ and
$v_i ^+\in E_i ^+$. Since the sets $\O_i$'s are disjoint, we
derive
$$
\sum_{i=1}^p\int_{\O_i}\mid \D v_i\mid^2\leq \sum_{i=1}^p\int_{\O_i}\mid \D
v\mid^2
$$
Thus
 \begin{align}\label{Q1}
 Q(v,v) & =\int_{R^n}\mid \D v\mid^2-\frac{n+4}{n-4}\sum_{i=1}^p
 \int_{R^n}\d_i ^{\frac{8}{n-4}} v^2\notag\\
 & =\int_{(\cup\O_i)^c}\mid\D v\mid^2+\sum_{i=1}^p\biggl(\int_{\O_i}\mid \D
 v\mid^2-\int_{\O_i}\mid \D v_i\mid^2\biggr)\notag \\
  & +\sum_{i=1}^p\biggl(\int_{\O_i}\mid \D v_i\mid^2
 -\frac{n+4}{n-4} \int_{R^n}\d_i ^{8/(n-4)} v_i ^2\biggr)-\frac{n+4}{n-4}\sum_{i=1}^p
 \int_{R^n}\d_i ^{\frac{8}{n-4}} (v^2-v_i ^2)
 \end{align}
 Observe that, using Lemmas \ref{54} and \ref{55}, we have
 \begin{align*}
 \int_{R^n} \mid \D v_i\mid^2 & -\frac{n+4}{n-4}\int_{R^n}\d_i^{\frac{8}{n-4}} v_i ^2=  \int_{R^n}\mid \D v_i ^+\mid^2+\int_{R^n}\mid \D v_i ^-\mid^2
 -\frac{n+4}{n-4}\int_{R^n}\d_i^{\frac{8}{n-4}} (v_i ^+)^2\\
 &
-\frac{n+4}{n-4}\int_{R^n}\d_i^{\frac{8}{n-4}} ((v_i ^-)^2+2v_i ^+v_i ^-)
   \geq \a_1\int_{R^n}\mid \D v_i ^+\mid^2+o(\int_{R^n}\mid \D
  v\mid^2)\\
 & \geq \frac{\a_1}{2}\int_{\O_i}\mid \D v_i \mid^2.
 \end{align*}
  We also have, using Lemma \ref{53},
 \begin{align}\label{v2}
 \int_{R^n}\d_i ^{\frac{8}{n-4}}(v^2-v_i
 ^2) & =\biggl[\int\mid v+v_i\mid^{\frac{2n}{n-4}}\biggr]^{\frac{n-4}{2n}}\biggl[\int\d_i
 ^{\frac{n+4}{n-4}}\mid v-v_i\mid\biggr]^{\frac{8}{n+4}}\biggl[\int\mid v-v_i
 \mid^{\frac{2n}{n-4}}\biggr]^{\frac{(n-4)^2}{2n(n+4)}}\notag\\
   & =o(\int\mid\D v\mid^2).
  \end{align}
   Thus, \eqref{Q1} becomes
  \begin{align}\label{Q2}
  Q(v,v) & \geq \int_{(\cup\O_i)^c}\mid \D
  v\mid^2+\sum_{i=1}^p\biggl[\int_{\O_i}\mid \D
  v\mid^2-\int_{\O_i}\mid \D
  v_i\mid^2\biggr]+\frac{\a_1}{2}\sum_{i=1}^p\int_{\O_i}\mid \D
  v_i\mid^2\notag \\
 & +o(\int_{R^n}\mid \D v\mid^2)
  \geq \int_{(\cup\O_i)^c}\mid \D
  v\mid^2+\frac{\a_1}{2}\sum_{i=1}^p\int_{\O_i}\mid \D
  v\mid^2+o(\int_{R^n}\mid \D v\mid^2)\notag\\
   & \geq \ov{\a}_0\int_{R^n}\mid \D v\mid^2.
  \end{align}
 Thus the proof of Proposition \ref{Q} is completed under Lemmas
 \ref{53}, \ref{54} and \ref{55}.
 \end{pf}

Next, we will come to the proofs of Lemmas \ref{53}, \ref{54} and \ref{55}.\\ \begin{pfn}{ Lemma \ref{54}}
Observe that the family of functions $\d_{(a,\l)}$ are the
solutions of the Yamabe problem on $\R^n$, that is, the functional
$$I(u)=\frac{1}{2}\int_{R^n}\mid \D u\mid^2-\frac{n-4}{2n}
\int_{R^n}\mid u\mid ^{\frac{2n}{n-4}}$$ has only the family of
functions ${\d}_{(a,\l)}$ as critical points. Those critical
points are degenerated and of index 1. The nullity space
 is of dimension $n+1$ and it is generated by the derivative of $\d_{(a,\l)}$ with respect to $\l$ and $a$.
 Furthermore, the set of negativity is generated by the function
 $\d:=\d_{(a,\l)}$. Let
$$F=span\{{\d},\frac{\partial {\d}}{\partial\l}, \frac{\partial
{\d}}{\partial (a)_i},i=1,...,n\}.
$$
Thus, on the orthogonal of $F$,
the second derivation of the functional $I$ on the point ${\d}$ is
positive definite. Therefore
\begin{eqnarray}\label{a1}
\exists \, \, \a_1>0 \, \, s.t.\, \, \forall\, v\in F^{\bot}\,
\mbox{we have } \, \mid\mid v\mid\mid^2-\frac{n+4}{n-4}\int_{S^n}
{\d}^{\frac{8}{n-4}}v^2\geq \a_1 \mid\mid v\mid\mid^2.
\end{eqnarray}
\end{pfn}\\
 \begin{pfn}{Lemma \ref{55}}
  We have
 $$
v_i ^-=a\d_i+b\l_i\frac{\partial \d_i}{\partial\l_i}+\sum_{j=1}^nc_j
 \frac{1}{\l_i}\frac{\partial \d_i}{\partial (a_i)_j}.
$$
 Multiply $v_i ^-$ by $\d_i ^{(n+4)/(n-4)}$ and integrate on
 $R^n$, using Lemma \ref{53}, we derive that
 \begin{eqnarray*}
 a\int_{R^n}\d_i ^{\frac{2n}{n-4}}=\int_{R^n}\d_i
 ^{\frac{n+4}{n-4}}v_i ^-=\int_{R^n}\d_i
 ^{\frac{n+4}{n-4}}(v_i ^--v)=O\biggl(
 (\sum\e_{ij}^{1/2})(\int_{R^n}\mid\D v\mid^2)^{1/2}\biggr).
 \end{eqnarray*}
 In the same way, we have
 \begin{align*}
 b\int_{R^n}\d_i ^{\frac{8}{n-4}}\mid \l_i\frac{\partial\d_i}
 {\partial\l_i}\mid^2 & =\int_{R^n}\d_i ^{\frac{8}{n-4}}\l_i\frac{\partial\d_i}
 {\partial\l_i}v_i ^-=\int_{R^n}\d_i ^{\frac{8}{n-4}}\l_i\frac{\partial\d_i}
 {\partial\l_i}(v_i ^--v)\\
  & =O\biggl(\int_{R^n}\d_i ^{\frac{n+4}{n-4}}\mid v_i -v\mid\biggr)
=O\biggl((\sum\e_{ij}^{1/2})(\int_{R^n}\mid\D
  v\mid^2)^{1/2}\biggr)
 \end{align*}
 \begin{align*}
 c_j\int_{R^n}\d_i ^{\frac{8}{n-4}}\mid \frac{1}{\l_i}\frac{\partial\d_i}
 {\partial (a_i)_j}\mid^2 & =\int_{R^n}\d_i ^{\frac{8}{n-4}}\frac{1}{\l_i}
 \frac{\partial\d_i}{\partial (a_i)_j}v_i ^-
 =\int_{R^n}\d_i ^{\frac{8}{n-4}}\frac{1}{\l_i}\frac{\partial\d_i}
 {\partial (a_i)_j}(v_i ^--v)\\
  & =O\biggl(\int_{R^n}\d_i ^{\frac{n+4}{n-4}}\mid v_i -v\mid\biggr)
=O\biggl((\sum\e_{ij}^{1/2})(\int_{R^n}\mid\D
  v\mid^2)^{1/2}\biggr).
 \end{align*}
 Using the fact that $\d_i$, $\l_i(\partial\d_i)/(\partial\l_i)$ and
 $\l_i ^{-1}(\partial\d_i)/(\partial (a_i)_j)$ have a constant norme, the lemma follows.
 \end{pfn}\\
  Before giving the proof of Lemma \ref{53}, we need the following
  lemma.
  \begin{lem}\label{56}
 Let $w\in H$ and $h$ satisfies
 $$\D^2 h =0\, \, on\, \, B_\l\, ,\, \D h=\D w \, \,
 and\, \, h=w\, \, in\, \, \partial B_\l,$$
 where $B_\l=\{x/\mid x\mid <\l\}$. We have
 \begin{eqnarray*}
 \int_{ B_\l}\d_{(0,1)}^{(n+4)/(n-4)}\mid h\mid \leq
 \frac{c}{\l^{(n-4)/2}}(\int_{R^n}\mid \D w\mid^2)^{1/2}.
 \end{eqnarray*}
 \end{lem}
 \begin{pf}
 First, observe that we have, for $w\in H$, if a function $u$ satisfies
 \begin{eqnarray}\label{36}
 \D^2 u =0\, \, on\, \, B_1\, ,\, \D u=\D w \, \,
 and\, \, u=w\, \, in\, \, \partial B_1,
 \end{eqnarray}
 thus
 \begin{eqnarray}\label{44}
 \int_{\partial B_1}\mid \D u\mid +\int_{\partial B_1}\mid
 u\mid\leq C (\int_{R^n}\mid\D w\mid^2)^{1/2}.
 \end{eqnarray}
 Notice that we can assume that $w\geq 0$ and $\D w\geq 0$. Indeed:
 \begin{rem}\label{rem}
 If the function $w$ does not satisfy "$-\D w\geq
 0$" and "$w \geq 0$", we can introduce the function
 $w'$ defined by
 $$w'\in H^1(\R^n) \mbox{ and } -\D w'=\mid\D w\mid \qquad on \quad R^n.$$
 Thus $w'\in H$ and it satisfies "$-\D w'\geq 0$" and
 "$w'\geq 0$" on $R^n$ and it easy to see that
 $$w'-w\geq 0,\quad , \quad \mid w-
 P w\mid \leq w'-P w'\quad and \quad \int_{R^n}\mid
 \D w'\mid^2=\int_{R^n}\mid \D w\mid^2, $$
 where $P$ is the projection operator on any subset.
 \end{rem}
 Hence if the lemma holds with $w\geq 0$ and $\D w\geq 0$, it will hold for all
 $w$. If we assume $w$ and $-\D w$ to be positive, the function $h$ will also
 be  positive. Then
 \begin{eqnarray}\label{41}
 \int_{B_\l}\d^{\frac{n+4}{n-4}}h =
 \int_{B_\l}\D^2 (\d-\th) h
  = \int_{\partial B_\l}\frac{\partial}{\partial n} (\D(\d-\th))h+
 \int_{\partial B_\l}\frac{\partial}{\partial n} (\d-\th)\D h
 \end{eqnarray}
 where $\th$ satisfies
 $$\D^2 \th =0\, \, on\, \, B_\l\, ,\, \D \th=\D \d \, \,
 and\, \, \th=\d\, \, in\, \, \partial B_\l$$
 It is easy to see that $\th$ is equal to
 $$\th=\frac{1}{2n}c_\l(\mid
 x\mid^2-\l^2)+\frac{c_0}{(1+\l^2)^{(n-4)/2}},$$
 where $c_0$ is defined in the definition of $\d$ and $c_\l$ is
 equal to
 $$c_\l=\D\d\mid_{\partial
 B_\l}=\frac{(n-4)c_0}{(1+\l^2)^{n/2}}(-n-2\l^2)\sim\frac{-c}{\l^{n-2}}$$
 (for $\l$ large). Therefore, we have
 \begin{eqnarray}\label{n1}
 \frac{\partial}{\partial n} (\d-\th)\mid_{\partial
 B_\l}=-\frac{(n-4)c_0\l}{(1+\l^2)^{(n-2)/2}}-\frac{\l}{n}c_\l=
 -\frac{(n-2)(n-4)c_0\l^3}{n(1+\l^2)^{n/2}}\sim -\frac{c}{\l^{n-3}}
 \end{eqnarray}
 and
 \begin{eqnarray}\label{n2}
 \frac{\partial}{\partial n} (\D(\d-\th))\mid_{\partial B_\l}=
 \frac{c_0(n-2)(n-4)\l}{(1+\l^2)^{\frac{n+2}{2}}}(n+2+2\l^2)\sim
 \frac{c}{\l^{n-1}}
 \end{eqnarray}
 (for $\l$ large). Using \eqref{n1} and \eqref{n2}, \eqref{41}
 becomes
 \begin{eqnarray}\label{11}
 \int_{B_\l}\d^{\frac{n+4}{n-4}}h \leq
 \frac{c}{\l^{n-1}}\int_{\partial B_\l}h +
 \frac{c}{\l^{n-3}}\int_{\partial B_{\l}}-\D h
 \end{eqnarray}
 Let $\ov{h}(x)=\l^{(n-4)/2}h(\l x)$ and $\ov{w}(x)=\l^{(n-4)/2}
 w(\l x)$. The function $\ov{h}$ satisfies \eqref{36} with $\ov{w}$
 instead of $w$. Thus, it satisfies \eqref{44} with $\ov{w}$.
 Observe that
 \begin{eqnarray}\label{45}
 \int_{\partial B_1}\ov{h}=\frac{1}{\l^{\frac{n+2}{2}}}\int_{\partial
 B_{\l}}h\,  ,\quad \int_{\partial B_1}\D\ov{h}=
 \frac{1}{\l^{\frac{n-2}{2}}}\int_{\partial B_{\l}}\D h\,  ,
 \quad \int_{R^n}\mid\D\ov{w}\mid^2=\int_{R^n}\mid\D w\mid^2.
 \end{eqnarray}
 Thus, using \eqref{n1}, \eqref{n2}, \eqref{11} and \eqref{45},
 the lemma follows.
 \end{pf}\\
 \begin{pfn}{Lemma \ref{53}}
 First, we assume that we have only two masses.
 Take $i=1$ in the Lemma. We can make a translation and a
 dilatation so that $\ov{\l}_1=1$ and $\ov{a}_1=0$. Let
 $$\ov{\varphi}=\l_1^{\frac{n-4}{2}}\varphi(\l_1x+a_1)\quad ,\quad
 \ov{\d}_j=\l_1^{\frac{n-4}{2}}\d_j(\l_1x+a_1)$$
 Notice that
 $$\ov{\e}_{12}=\e_{12}\quad ,\quad \ov{\l}_1=1\quad ,\quad
 \ov{\l}_2=\frac{\l_2}{\l_1}\quad , \quad \ov{a}_1=0\quad , \quad
 \ov{a}_2=\l_1(a_2-a_1)$$
 Assume first that $\l_1\geq \l_2$, hence $\ov{\l}_2\leq 1$. Then
 $$\ov{\O}_1=\{x/\mid x\mid<(8\e_{12}^{1/(n-4)})^{-1}\}$$
 Let $\ov{\var}_1=\ov{\var}-\ov{h}$ with $\D^2 \ov{h}=0$ in
 $\ov{\O}_1$, $\D\ov{h}=\D\ov{\var}$ and $\ov{h}=\ov{\var}$
 on $\partial \ov{\O}_1$, we have
 \begin{align}\label{52}
 \int_{R^n}\ov{\d}_1^{\frac{n+4}{n-4}}\mid\ov{\var}-\ov{\var}_1\mid
  & \leq \int_{\ov{\O}_1}+\int_{\mid
 x\mid\geq(8\e_{12}^{1/(n-4)})^{-1}}\notag \\
 & \leq c\e_{12}^{1/2}(\int_{R^n}\mid\D \ov{\var}\mid^2)^{1/2}+
 \e_{12}^{\frac{n+4}{2(n-4)}}(\int_{R^n}\mid\D
 \ov{\var}\mid^2)^{1/2}\notag \\
  & \leq c(\e_{12}^{1/2}+ \e_{12}^{\frac{n+4}{2(n-4)}})
 (\int_{R^n}\mid\D {\var}\mid^2)^{1/2}
 \end{align}
(using Lemma \ref{56} and Holder's inequality.) Observe that
$$
\int_{R^n}\ov{\d}_1^{\frac{n+4}{n-4}}\mid\ov{\var}-\ov{\var}_1\mid
=\int_{R^n}{\d}_1^{\frac{n+4}{n-4}}\mid{\var}-{\var}_1\mid.
$$
 Thus, the proof is completed in this case (the case where
 $\l_2\leq \l_1$). We will now see the other case i.e.
 $\l_1\leq \l_2$. Thus
$$
\ov{\O}_1=\{x/\mid x\mid<(8\e_{12}^{1/(n-4)})^{-1}\, \mbox{ and }\,
\mid x-\ov{a}_2\mid>(8\ov{\l}_2\e_{12}^{1/(n-4)})^{-1}\}.
$$
 Let
 $$
\wtilde{\O}=\{x/\mid x\mid<(8\e_{12}^{1/(n-4)})^{-1}\}\quad
 ,\quad \wtilde{W}=\{x/\mid
 x-\ov{a}_2\mid>(8\ov{\l}_2\e_{12}^{1/(n-4)})^{-1}\}.
$$
 Observe that we have $\partial\ov{\O}_1\subset\partial\wtilde{\O}
 \cup \partial\wtilde{W}$. We define $\wtilde{\var}_1$ to be the
 projection of $\ov{\var}$ on $\wtilde{\O}$ and $\wtilde{\psi}_1$ to
 be the projection of $\ov{\var}$ on $\wtilde{W}$.\\
 In the following, we will assume that $-\D\ov{\var}\geq 0$ and
 $\ov{\var}\geq 0$. The general case can be deduced by
 Remark \ref{rem}. \\
 Hence we derive
 \begin{eqnarray}\label{61'}
 \mid \ov{\var}-\ov{\var}_1\mid \leq (\ov{\var}-\wtilde{\var}_1)
 +(\ov{\var}-\wtilde{\psi}_1)\qquad in \quad R^n
 \end{eqnarray}
 and thus
 \begin{eqnarray}\label{62}
 \int_{R^n}\ov{\d}_1^{\frac{n+4}{n-4}}\mid \ov{\var}-\ov{\var}_1\mid \leq
  \int_{\wtilde{\O}}\ov{\d}_1^{\frac{n+4}{n-4}}(\ov{\var}-\wtilde{\var}_1)
 +\int_{\wtilde{W}}\ov{\d}_1^{\frac{n+4}{n-4}}(\ov{\var}-\wtilde{\psi}_1)
 +\int_{\wtilde{\O}^c}\ov{\d}_1^{\frac{n+4}{n-4}}\ov{\var}
 +\int_{\wtilde{W}^c}\ov{\d}_1^{\frac{n+4}{n-4}}\ov{\var}.
 \end{eqnarray}
 As in \eqref{52}, using Holder's inequality, we have
 $$\int_{\wtilde{\O}^c}\ov{\d}_1^{\frac{n+4}{n-4}}\ov{\var}\leq
 c\e_{12}^{\frac{n+4}{2(n-4)}}(\int_{R^n}\mid\D\ov{\var}\mid^2)^{1/2}.$$
 We estimate now
 $\int_{\wtilde{W}^c}\ov{\d}_1^{\frac{n+4}{n-4}}\ov{\var}$.
 As in \cite{B1}, we prove that
 \begin{eqnarray}\label{64}
 \ov{\d}_1=\frac{c_0}{(1+\mid x\mid^2)^{\frac{(n-4)}{2}}}\leq
 c_0\frac{2^{\frac{n-4}{2}}\ov{\l}_2^{\frac{n-4}{2}}}{(1+\ov{\l}_2^2\mid
 x-\ov{a}_2\mid^2)^{\frac{n-4}{2}}} \qquad \mbox{ if }\qquad \ov{\l}_2\mid
 x-\ov{a}_2\mid\leq \frac{1}{8\e_{12}^{\frac{1}{n-4}}}.
\end{eqnarray}
 We then have
 \begin{align}\label{70}
 \int_{\wtilde{W}^c}\ov{\d}_1^{\frac{n+4}{n-4}}\ov{\var} & \leq c
 \biggl(\int_{\wtilde{W}^c} \ov{\d}_1^{\frac{2n}{n-4}}\biggr)^{\frac{n+4}{2n}}
 \biggl(\int_{R^n} \mid\D\ov{\var}\mid^2\biggr)^{\frac{1}{2}}
  \leq \biggl(\int_{\wtilde{W}^c}(\ov{\d}_1\ov{\d}_2)^{\frac{n}{n-4}}
 \biggr)^{\frac{n+4}{2n}}\biggl(\int_{R^n}
 \mid\D\ov{\var}\mid^2\biggr)^{\frac{1}{2}}\notag \\
 & \leq c \e_{12}^{\frac{(n+4)(n-2)}{2n(n-4)}}\biggl(\int_{R^n}
 \mid\D\ov{\var}\mid^2\biggr)^{\frac{1}{2}} \leq c\e_{12}^{\frac{1}{2}}
 \biggl(\int_{R^n}\mid\D\ov{\var}\mid^2\biggr)^{\frac{1}{2}}.
 \end{align}
 Using Lemma \ref{56}, and as in \eqref{52}, we have
 \begin{eqnarray}\label{72}
  \int_{\wtilde{\O}}\ov{\d}_1^{\frac{n+4}{n-4}}(\ov{\var}-\wtilde{\var}_1)
  \leq c\e_{12}^{\frac{1}{2}}\biggl(\int_{R^n}
 \mid\D\ov{\var}\mid^2\biggr)^{\frac{1}{2}}.
 \end{eqnarray}
 It remains to estimate
 $\int_{\wtilde{W}}\ov{\d}_1^{\frac{n+4}{n-4}}(\ov{\var}-\wtilde{\psi}_1)$
 \begin{align}\label{e:618}
 \int_{\wtilde{W}}\ov{\d}_1^{\frac{n+4}{n-4}}(\ov{\var}-\wtilde{\psi}_1) & =
 \int_{\wtilde{W}}\D^2\ov{\d}_1(\ov{\var}-\wtilde{\psi}_1)=
 \int_{\partial \wtilde{W}}\frac{\partial}{\partial
 n}(\ov{\d}_1-\ov{\th}_1)\D\ov{\var}+\int_{\partial
 \wtilde{W}}\frac{\partial}{\partial
 n}(\D(\ov{\d}_1-\ov{\th}_1)]\ov{\var}\notag \\
 & \leq \sup_{\partial\wtilde{W}}\mid\frac{\partial}{\partial
 n}(\ov{\d}_1-\ov{\th}_1)\mid \int_{\partial\wtilde{W}}\D\ov{\var}
 +\sup_{\partial\wtilde{W}}\mid\frac{\partial}{\partial
 n}(\D(\ov{\d}_1-\ov{\th}_1))\mid \int_{\partial\wtilde{W}}\ov{\var},
 \end{align}
 where $\ov{\th}_1$ is the projection of $\ov{\d}_1$ on
 $\wtilde{W}$. Now, we need to estimate the normal derivatives which appear in \eqref{e:618}. For this effect, let us introduce the Green's function $G_{\wtilde W}$ which satisfies
$$
\D^2 G_{\wtilde W}(x,.)=\d_x , \mbox{ in } \wtilde W, \quad \D G_{\wtilde W}=G_{\wtilde W}=0 \mbox{ on } \partial \wtilde W.
$$
Thus for any function $u$ we have
$$
u(y)=\int_{\wtilde W}G_{\wtilde W}\D^2 u +\int_{\partial \wtilde W}\frac{\partial}{\partial \nu}G_{\wtilde W}\D u +\int_{\partial \wtilde W}\frac{\partial}{\partial \nu}(\D G_{\wtilde W}) u.
$$
Observe that $(\ov{\d}_1-\ov{\theta}_1)$ satisfies
$$
\D^2 (\ov{\d}_1-\ov{\theta}_1)=\ov{\d}_1^{\frac{n+4}{n-4}} \mbox{ in }\wtilde W, \quad \D(\ov{\d}_1-\ov{\theta}_1)=\ov{\d}_1-\ov{\theta}_1=0 \mbox{ on } \partial\wtilde W.
$$
Thus, we derive
\begin{eqnarray}\label{e:*}
(\ov{\d}_1-\ov{\theta}_1)(y)=\int_{\wtilde W}G_{\wtilde W}(x,y)\ov{\d}_1^{\frac{n+4}{n-4}}(x),\,  |\frac{\partial}{\partial \nu_y}(\ov{\d}_1-\ov{\theta}_1)(y)|=\int_{\wtilde W}\frac{\partial}{\partial\nu_y}G_{\wtilde W}(x,y)\ov{\d}_1^{\frac{n+4}{n-4}}(x).
\end{eqnarray}
But we have
$$
G_{B^c(o,1)}(x,y)=\frac{1}{|x|^{n-4}}G_{B(o,1)}(\frac{x}{|x|^2},\frac{y}{|y|^2}),$$
$$G_{\wtilde W}(x,y)=\frac{1}{r^{n-4}}G_{B^c(o,1)}(\frac{x}{r},\frac{y}{r})=\frac{1}{|x|^{n-4}}G_{B(o,1)}(\frac{rx}{|x|^2},\frac{ry}{|y|^2}).
$$
Let $y\in \partial\wtilde W$ and let $\pi_y$ be the half space which contains $B(o,1)$ and satisfies $y\in \partial \pi_y$, then we have
$$
G_{B(o,1)}(x,y)\leq G_{\pi_y}(x,y) \mbox{ and } |\frac{\partial}{\partial \nu_y} G_{B(o,1)}(x,y)|\leq |\frac{\partial}{\partial\nu_y} G_{\pi_y}(x,y)| \leq \frac{c}{|x-y|^{n-3}}
$$
and therefore, since $y\in \partial\wtilde W$,
\begin{align*}
 |\frac{\partial}{\partial \nu_y} G_{\wtilde W}(x,y)| & \leq \frac{c}{r|x|^{n-4}}\frac{1}{|\frac{rx}{|x|^2}-\frac{y}{r}|^{n-3}}\leq \frac{c|x|}{r|x-y|^{n-3}}\\
 & \leq \frac{c}{r|x-y|^{n-4}}+\frac{c}{|x-y|^{n-3}}.
\end{align*}
Using \eqref{e:*} we derive
\begin{align*}
|\frac{\partial}{\partial \nu_y}(\ov{\d}_1-\ov{\theta}_1)(y)| & \leq \int_{R^n}\frac{c\d^{\frac{n+4}{n-4}}}{r|x-y|^{n-4}}+c\int_{R^n}\frac{\d^{\frac{n+4}{n-4}}(x-y)}{|x|^{n-3}}\\
 & \leq \frac{c}{r}\d(y)+\frac{4c}{(1+r^2)^{1/2}}\int_{|x|^2\geq (1+r^2)/4}\frac{\d^{\frac{n+4}{n-4}}(x-y)}{|x|^{n-4}}+c\int_{4|x|^2\leq (1+r^2)}\frac{\d^{\frac{n+4}{n-4}}(y)}{|x|^{n-3}}\\
 & \leq \frac{c}{r}\d(y)+c\d(y)^{\frac{n+4}{n-4}}(1+r^2)^{3/2} \leq \frac{c}{r(1+r^2)^{(n-4)/2}}.
\end{align*}
For the second term, we introduce  the Green's function  $\wtilde G_{\wtilde W}$ for $-\D$, i.e. $\wtilde G_{\wtilde W}$ satisfies
$$
-\D \wtilde G_{\wtilde W}(x,.)=\d_x \mbox{ in } \wtilde W, \quad \wtilde G_{\wtilde W}=0 \mbox{ on } \partial \wtilde W.
$$
By the same argument we prove that
$$\bigg|\frac{\partial}{\partial \nu_y}\wtilde G_{\wtilde W}\bigg|\leq \frac{c}{r|x-y|^{n-2}}+\frac{c}{|x-y|^{n-1}}.
$$
Arguing as above, for $g=\D(\ov{\d}_1-\ov{\theta}_1)$, we have
\begin{align*}
\bigg|\frac{\partial g(y)}{\partial \nu_y}\bigg| & =\bigg|\frac{\partial}{\partial \nu_y}(\D(\ov{\d}_1-\ov{\theta}_1)(y)\bigg|\leq \int_{R^n}\frac{c\d^{\frac{n+4}{n-4}}}{r|x-y|^{n-2}}+\int_{R^n}\frac{\d(x-y)^{\frac{n+4}{n-4}}}{|x|^{n-1}}\\
 & \leq \frac{c}{r}g(y)+\frac{c}{(1+r^2)^{1/2}}g(y)+\d(y)^{\frac{n+4}{n-4}}(1+r^2)^{1/2}\leq \frac{c}{r(1+r^2)^{\frac{n-2}{2}}}.
\end{align*}
Thus \eqref{e:618} becomes
$$
 \int_{\wtilde{W}}\ov{\d}_1^{\frac{n+4}{n-4}}(\ov{\var}-\wtilde{\psi}_1) \leq
\frac{c}{r(1+r^2)^{\frac{n-4}{2}}}\biggl(\int_{\partial \wtilde W}\D\ov{\var}+\frac{1}{r^2}\int_{\partial \wtilde W}\ov{\var}\biggr).
$$
Using \eqref{44} and \eqref{45}, we derive
$$
\int_{\wtilde{W}}\ov{\d}_1^{\frac{n+4}{n-4}}(\ov{\var}-\wtilde{\psi}_1) \leq \frac{c}{r(1+r^2)^{\frac{n-4}{2}}}r^{\frac{n-2}{2}}\int_{R^n}|\D \var|^2\leq \frac{cr^{\frac{n-4}{2}}}{(1+r^2)^{\frac{n-4}{2}}}||\var||^2.
$$
Recall that $r=(8\ov{\l}_2\e_{12}^{1/(n-4)})^{-1}$. If $|\ov{x}_2|<1$ then $\e_{12}\sim \ov{\l}_2^{(4-n)/2}$ and therefore $r^{(n-4)/2}\leq \ov{\l}_2^{(4-n)/4}\leq c \e_{12}^{1/2}$. In the other case, that is,  $|\ov{x}_2|\geq 1$, we have $\e_{12}^{-1}\sim (\ov{\l}_2|\ov{x}_2|^2)^{(n-4)/2}$ and therefore $r^{(4-n)/2}\leq (\ov{\l}_2/|\ov{x}_2|^2)^{(n-4)/2}\leq \e_{12}^{1/2}$. Thus, in all cases we obtain
$$
 \int_{\wtilde{W}}\ov{\d}_1^{\frac{n+4}{n-4}}(\ov{\var}-\wtilde{\psi}_1) \leq \e_{12}^{1/2}||\var||^2.
$$
 This completes the proof in the case where we are dealing with
 two points.\\
 In the general case, one introduces the sets, assuming
 $\l_1=1$, $a_1=O$ and $\var\geq 0$
 $$W_i=\{x\in R^n/\mid x\mid<\e_{1i}^{-1/(n-4)},\quad \mid
 x-a_i\mid>\l_i ^{-1}\e_{1i}^{-1/(n-4)}\, if\, \l_i>1\}$$
 Then $\partial\O_1\subset\cup\partial W_i$.
 Let $\var_1$ be the projection of $\var$ on $\O_1$ and
 $\wtilde{\var}_i$ be the projection of $\var$ on $W_i$. Then the
 above arguments, in the case of two points, imply
 \begin{eqnarray}\label{86}
 \int_{W_i}\d_1^{\frac{n+4}{n-4}}(\var-\wtilde{\var}_i)\leq
 c\e_{1i}^{1/2}(\int_{R^n}\mid\D\var\mid^2)^{1/2}
 \end{eqnarray}
 \begin{eqnarray}\label{87}
 \int_{W_i ^c}\d_1^{\frac{n+4}{n-4}}\var\leq
 c\e_{1i}^{1/2}(\int_{R^n}\mid\D\var\mid^2)^{1/2}
 \end{eqnarray}
 From \eqref{86} and \eqref{87} the general case follows.
 \end{pfn}
 \subsection{ Some estimates }
 \mbox{}
In this subsection, we collect some technical estimates of the different integral quantities which occur in the paper. The proof of these estimates are similar  to their analogous for Laplacian in \cite{B1} and \cite{R}.
 \begin{lem}\label{lem1}
Let $a\in S^n$ and $\l > 0$ large enough. Using the stereographic projection $\pi_{-a}$ the function $\wtilde\d_{(a,\l)}$ will be transformed to $\d_{(0,\l)}$ (see \cite{BB}). Furthermore, we have
$$
\int_{S^n}L\wtilde{\d}.\wtilde{\d} =
\int_{S^n}\wtilde{\d}^{\frac{2n}{n-4}}=\int_{R^n}\d
^{\frac{2n}{n-4}}=S_n, \quad <\wtilde{\d},\l \frac{\partial \wtilde{\d}}{\partial \l}>=0, \quad <\wtilde{\d},\frac{1}{\l} \frac{\partial \wtilde{\d}}{\partial  a}>=0.
 $$
\end{lem}
\begin{lem}\label{lem2}
For $a_1$, $a_2\in S^n$, $\l_1$, $\l_2 > 0$ large enough, let $b\in S^n$ such that $d(a_1,b)=d(a_2,b)$. Using the stereographic projection $\pi_{-b}$, the function $\wtilde\d_{(a_i,\l_i)}$ will be transformed to $\d_{(\tilde {a}_i,\tilde{\l}_i)}$ with
$$\tilde{a}_i=\frac{(\l_i ^2-1)Proj_{R^n}a_i}{2+(\l_i ^2-1)(1-\cos (\theta_0))}, \quad \tilde{\l}_i=\frac{2+(\l_i ^2-1)(1-\cos (\theta_0))}{2\l_i}, \quad \theta_0=\pi-d(a_i,b)
$$
(see \cite{BB}). Furtheremore, we have for $i\ne j$,
$$\int_{S^n}L\wtilde{\d}_i.\wtilde{\d}_j=\int_{S^n}\wtilde{\d}_i
^{\frac{n+4}{n-4}}\wtilde{\d}_j=\int_{R^n}\d_i
^{\frac{n+4}{n-4}}\d_j=c_1\tilde\e _{ij}+o(\tilde\e_{ij})=c_1\e _{ij}+o(\e_{ij})$$
 where
$$\e_{ij}=\biggl(\frac{\l_i}{\l_j}+\frac{\l_j}{\l_i}+\frac{\l_i\l_j}{2}(1-\cos d(a_i,a_j))\biggr)^{-\frac{n-4}{2}} \mbox { and } \tilde\e_{ij}=\biggl(\frac{\tilde\l_i}{\tilde\l_j}+\frac{\tilde\l_j}{\tilde\l_i}+\tilde\l_i\tilde\l_j|\tilde a_i-\tilde a_j|^2\biggr)^{-\frac{n-4}{2}}
$$
and  $c_1=\b_n^{2n/(n-4)}\int_{R^n}(1+\mid
 x\mid^2)^{-(n+4)/2}$. If $n=6$, $c_1=\frac{\b_6^{6}}{24}w_5$, $\b_n$ is defined in the definition
 of $\wtilde{\d}$ and $w_5$ is the volume of the five dimensional sphere.
\end{lem}
\begin{lem}\label{lem3}
We have the following estimates.
\begin{align*}
\int_{S^n}K\wtilde{\d}^{\frac{2n}{n-4}}&=\int_{R^n}\wtilde{K}\d
^{\frac{2n}{n-4}}=K(a)S_n+c_2\frac{4\D K(a)}{\l ^2}+O(\frac{1}{\l
^3})\\
\int_{S^n}K \wtilde{\d}^{\frac{n+4}{n-4}}\l \frac{\partial
\wtilde{\d}}{\partial \l}&=-\frac{n-4}{n}c_2\frac{4\D
K(a)}{\l^2}+O(\frac{1}{\l^3})\\
\int_{S^n}K \wtilde{\d}^{\frac{n+4}{n-4}}\frac{1}{\l} \frac{\partial
\wtilde{\d}}{\partial a}&=c_3\frac{\n K(a)}{\l}+O(\frac{1}{\l^2}),
\end{align*}
 where $c_2=\frac{1}{2n}\int_{R^n}\mid x\mid^2\d_{(O,1)}^{2n/(n-4)}$.
 If $n=6$, $c_2 =\frac{\b_6^{6}}{480}w_5$.
\end{lem}
\begin{lem}\label{lem4}
 For $i\ne j$, we have the following estimates
\begin{align*}
\int_{S^n}K\wtilde{\d}_i
^{\frac{n+4}{n-4}}\wtilde{\d}_j&=\int_{R^n}\wtilde{K}\d _i
^{\frac{n+4}{n-4}}\d _j=c_1K(a_i)\e_{ij}+o(\e_{ij}+\frac{1}{\l_i ^2})\\
\int_{S^n}(\wtilde{\d}_i\wtilde{\d}_j)^{\frac{n}{n-4}}&=O(\e
_{ij}^{\frac{n}{n-4}}log\e _{ij}^{-1})\\
<\wtilde{\d}_j,\l_i \frac{\partial \wtilde{\d}_i}{\partial
\l_i}> &= c_1 \l_i \frac{\partial \e_{ij}}{\partial \l_i}
+o(\e_{ij})\\
\int_{S^n}K \wtilde{\d}_j^{\frac{n+4}{n-4}}\l_i \frac{\partial
\wtilde{\d}_i}{\partial \l_i}&=c_1K(a_j)\l_i \frac{\partial
\e_{ij}}{\partial \l_i} +o(\e_{ij}+\frac{1}{\l_j^2})\\
\frac{n+4}{n-4}\int_{S^n}K \wtilde{\d}_i^{\frac{8}{n-4}}\l_i
\frac{\partial \wtilde{\d}_i}{\partial \l_i}\wtilde{\d}_j&=
c_1K(a_i)\l_i \frac{\partial \e_{ij}}{\partial \l_i}+
o(\e_{ij}+\frac{1}{\l_i ^2}).
\end{align*}
\end{lem}
{\bf Acknowledgements.} Part of this work was done when the authors enjoyed the hospitality of Rutgers Mathematics Department at New Brunswick (New Jersey, USA). They would like to thank the Mathematics Department for its warm hospitality. The authors also thank Professor Abbas Bahri for his encouragement and constant support over the years.

\end{document}